\documentclass[11pt,a4paper]{article} 
\usepackage{amsmath,pstricks,latexsym,theorem,epsfig} 
\newtheorem{Theorem}{Theorem}[section] 
\newtheorem{Definition}{Definition}[section] 
 
\newtheorem{Lemma}{Lemma}[section]

\theorembodyfont{\rmfamily} 
\newtheorem{Remark}{Remark}[section] 
\newtheorem{Example}{Example}[section]

\def \R{I\!\!R}

\def\11{1\!\!1}

\newcommand{\eps}       {\varepsilon} 
 
\newcommand{\ba}{\begin{array}} 
\newcommand{\ea}{\end{array}}

\begin{document} 

\title{ Uniqueness  results for convex Hamilton-Jacobi equations 
under  $p>1$ growth conditions on data} 

\author{Francesca Da Lio$^{(1)}$ \& 
Olivier Ley$^{(2)}$} 

\addtocounter{footnote}{1} 
\footnotetext{Dipartimento di Matematica Pura e Applicata, Via Trieste, 63, 
35121 Padova, Italy.} 

\addtocounter{footnote}{1} \footnotetext{Laboratoire de 
Math\'ematiques et Physique Th\'eorique (UMR CNRS 6083). F\'ed\'eration
Denis Poisson (FR~2964). 
Universit\'e Fran\c{c}ois Rabelais Tours. Parc de 
Grandmont, 37200 Tours, France.} 


\maketitle 

\begin{abstract}
Unbounded stochastic control problems may lead to
Hamilton-Jacobi-Bellman equations whose Hamiltonians
are not always defined, especially when the diffusion term
is unbounded with respect to the control.
We obtain existence and uniqueness of viscosity solutions
growing at most like $o(1+|x|^p)$ at infinity for such
HJB equations and more generally for degenerate parabolic
equations with a superlinear convex gradient nonlinearity.
If the corresponding control problem has a bounded diffusion
with respect to the control, then our results apply to
a larger class of solutions, namely those growing
like $O(1+|x|^p)$ at infinity. This latter case encompasses
some equations related to backward stochastic differential equations.
\end{abstract}

{\small
}

\section{Introduction} 

In the joint paper \cite{dll06} the authors obtain a comparison result 
between semicontinuous viscosity solutions, neither bounded
from below nor from above, growing  at 
most quadratically in the state variable, of second order degenerate parabolic equations of the form
\begin{equation}\label{hjbintr}
\left\{\begin{array}{ll}
\displaystyle\frac{\partial u}{\partial t}
+ H(x,t,Du, D^2 u) = 0 & \mbox{in  $\R^N\times (0,T),$}\\[3mm]
u (x,0)=\psi (x)& \mbox{in  $\R^N,$}\end{array}\right.
\end{equation}
where $N\geq 1,$ $T>0.$ The unknown $u$ is a real-valued function defined in $\R^N\times[0,T],$
$Du$ and $D^2u$ denote respectively its gradient and Hessian matrix
and $\psi$ is a given initial condition.
The Hamiltonian $H~\colon\R^N\times[0,T]\times \R^N\times{\cal{S}}_{N}(\R)\to\R$
has the form
\begin{eqnarray} \label{expression-H-intr}
H(x,t,q,X) = \mathop{\rm sup}_{\alpha\in A}
\left\{ -\langle b(x,t,\alpha),  q\rangle - \ell(x,t,\alpha)
-{\rm Trace}\left[ \sigma(x,t, \alpha)\sigma^T(x,t,\alpha)
 X \right]
\right\}.
\end{eqnarray}
Note that $H$ is convex with respect to $q.$
The key assumptions in the paper \cite{dll06} are that $A$ is an unbounded control set,
the functions $b$ and $\ell$ grow respectively at most linearly and quadratically 
with respect to both the control and the state.
Instead the functions $\sigma$ is assumed to grow at most linearly
with respect to the state and is bounded with respect to the control.
(In fact, in~\cite{dll06}, we consider more general equations of Isaacs
type by adding a concave Hamiltonian $G$ with bounded control, see
Remark~\ref{ajout-G}. To simplify the exposition we take $G\equiv 0$ here.)
\smallskip

In the present work, we extend the results of~\cite{dll06} in two
directions.
\smallskip

The first issue 
is to obtain a comparison result for unbounded solutions under the weaker assumption that 
the  diffusion matrix $\sigma$ is unbounded also with respect to the control.
The main difficulty is that the Hamiltonian $H$ 
may not be continuous. To illustrate this fact, 
consider for instance the case where $A=\R^N$, 
$b=\alpha,$ $\sigma =|\alpha|I$ and $\ell=|\alpha|^2$. The Hamiltonian $H$ becomes
\begin{eqnarray} \label{H-est-infini}
\sup_{\alpha\in \R^N}\{-\langle\alpha, Du\rangle
-|\alpha|^2-\frac{|\alpha|^2}{2}\Delta u\},
\end{eqnarray}
which is $+\infty$ as soon as $\Delta u< -2.$ 
This example is motivated by the well-known Stochastic Linear Quadratic problem, see for 
instance Bensoussan \cite{bensoussan82},
Fleming and Rishel \cite{fr75}, Fleming and Soner \cite{fs93}, {\O}ksendal \cite{oksendal98},
Yong and Zhou \cite{yz99} and the references therein for an overview of this problem.
The usual way to deal with such a problem is to plug into the equation
value functions $V$ of particular form (for instance quadratic in space)
for which one knows that $H(x,t,V,DV)$ is defined.
It leads to some ordinary
differential equations of Ricatti type which allow to identify
precisely the value function (see~\cite{yz99}).
Another way is to replace the Hamilton-Jacobi-Bellman equation
by a variational inequality, see Barles~\cite{barles90b} for instance.
Our aim is to study directly the PDE~\eqref{hjbintr}
without any a priori knowledge on the value function. Indeed, for
general datas, one does not expect explicit formula
for the value function.
\smallskip

We overcome the above difficulty in noticing that it is possible
to formulate the definition of viscosity solutions for HJB in a new
way without writing the  ``sup'' in~\eqref{expression-H-intr}, 
see Definition~\ref{newdef}. 
It provides a precise definition of solutions for~\eqref{hjbintr}
even in cases like~\eqref{H-est-infini}.  Let us stress that it is not
a new definition of viscosity solutions but only a new formulation.
Using this formulation,
we prove a comparison result for solutions in the class
of functions growing 
at most like $o(1+|x|^p)$ at infinity.
It provides new results for Stochastic Linear Quadratic type problems
(in this case, $p=2$) but, unfortunately, we are not able to treat
the classical Stochastic Linear Quadratic type problem with terminal
cost $\psi(x)=|x|^2$ since it requires a comparison in the class
$O(1+|x|^2).$ Nevertheless, our results apply
to very general datas (not only polynomials of degree 1 or 2
in $(x,\alpha)$), see Example~\ref{slqp}.
\smallskip

The second issue of our work is to extend the results of~\cite{dll06}
for $p$-growth type conditions on the datas and the solutions and
for more general equations with an additional nonlinearity $f$ which
is also convex with respect to the gradient and depends on $u.$ 
The motivation comes from PDEs arising in
the context of backward stochastic differential equations (BSDEs in short).
\smallskip

 In the framework of BSDEs, 
one generally considers forward-backward systems of the form
 \begin{equation}\label{diffuseq} 
\left\{\begin{array}{ll} 
 dX_s^{x,t}=b(X_s^{x,t},s)ds+\sigma(X^{t,x}_s,s)dW_s,& t\leq s\leq T,\\
 X_t^{t,x}=x,\end{array}\right.
 \end{equation}
  \begin{equation}\label{bsde} 
\left\{\begin{array}{ll} 
 -dY_s^{x,t}=f(X_s^{x,t},s,Y_s^{x,t},Z_s^{x,t})ds- Z_s^{x,t}dW_s,& t\leq s\leq T,\\
 Y_T^{x,t}=\psi(x),\end{array}\right.
 \end{equation}
where $(W_s)_{s\in[0,T]}$  is standard Brownian motion on a probability 
space $(\Omega,{\cal{F}},({\cal{F}}_t)_{t\in[0,T]},P)$, with $ ({\cal{F}}_t)_{t\in[0,T]}$ 
the standard Brownian filtration. (Note that $b$ and $\sigma$ do not
depend on the control). The diffusion 
\eqref{diffuseq} is associated with the second-order elliptic operator $L$  defined by
$$
Lu=-\frac{1}{2}{\rm Trace} (\sigma\sigma^TD^2u)-\langle b(x,t),Du\rangle.$$
The forward-backward 
system~\eqref{diffuseq}-\eqref{bsde}  is formally connected to the 
PDE
\begin{equation}\label{semihjb} 
\left\{\begin{array}{ll} 
-\displaystyle\frac{\partial u}{\partial t} 
+Lu-f(x,t,u,s(x,t)Du) 
  = 0 & \mbox{in  $\R^N\times (0,T)$}\\[3mm] 
u (x,T)=\psi (x) & \mbox{in  $\R^N.$}\end{array}\right.
\end{equation} 
by the nonlinear Feynman-Kac formula
\begin{eqnarray}\label{nlfk}
u(x,t)=Y_t^{x,t}  \quad \text{for all $(x,t)\in\R^N\times [0,T].$}
\end{eqnarray}
 We recall that nonlinear BSDEs with Lipschitz continuous  
coefficients were first introduced by Pardoux and Peng~\cite{pp90}, 
who proved existence and uniqueness.
Their results were extended by  Kobylanski~\cite{kobylanski00} 
for bounded solutions in the
case of coefficients $f$ having a quadratic growth in the gradient. 
Briand and Hu~\cite{bh08} generalized this latter result to the case
of solutions which are $O(1+|x|^p)$, as $|x|\to\infty,$  with $1\leq p<2.$ 
In all these works, the connection with viscosity solutions to 
the related PDE
\eqref{semihjb} is established: $u$ defined by~\eqref{nlfk}
is a viscosity solution of~\eqref{semihjb}.
\smallskip

Our aim is to prove the analytical counterpart of their
results. More precisely, we want to prove the existence and uniqueness
of the solution of~\eqref{semihjb} under the assumptions
of~\cite{bh08}.
\smallskip

Let us turn to a more precise description of our results.
We consider equations of the form
\begin{equation}\label{hjb} 
\left\{\begin{array}{ll} 
\displaystyle\frac{\partial u}{\partial t} 
+ H(x,t,Du, D^2 u) +f(x,t,u,s(x,t)Du) 
  = 0 & \mbox{in  $\R^N\times (0,T),$}\\[3mm] 
u (x,0)=\psi (x)& \mbox{in  $\R^N,$}\end{array}\right.
\end{equation} 
 where $H$ is given by \eqref{expression-H-intr} with $A$ unbounded,
$f\colon\R^N\times[0,T]\times\R^N\to\R$ is   
continuous and convex in the gradient and $s$ is bounded.
We look for solutions with $p>1$ growth assumptions 
(see~\eqref{classeC} and~\eqref{cprime}) and
both $H$ and $f$ satisfies some $p'$ growth assumptions,
where $p'=p/(p-1)$ is the conjugate of $p.$
See {\bf (A)}, {\bf (B)}, {\bf (C)} for the precise assumptions.  Let us
mention that the typical case we want to deal with is
$$
f(x,t,u,s(x,t)Du) = |s(x,t)Du|^{p'}, \quad p'>1,
$$
and the presence of $x$ in the power-$p'$ term is delicate
to treat (especially when doubling the variables in viscosity
type's proofs, see the proof of Lemma~\ref{linearisation}).
The $u$-dependence in $f$  means that $f(x,t,u,s(x,t)Du)$
may not be on the form of $H$ and induces some technical difficulties.
\smallskip

Section~\ref{sec:unbd} is devoted to the case with diffusion
matrices $\sigma$ which depend on the control in an unbounded way, 
see condition~\eqref{sig-non-borne}. The compensation
to  this condition with respect to~\cite{dll06}
(where $\sigma$ was assumed to be bounded with respect to
the control and $p=2$) is that 
we prove the comparison result Theorem~\ref{thm-unicite}
for semicontinuous sub- and supersolutions of~\eqref{hjb}
growing at most like $o(1+|x|^p)$ as
 $|x|\to\infty$ (instead of  $O(1+|x|^2)$ in~\cite{dll06}).
So far it remains an open question to know if there is uniqueness
 in the larger  class $O(1+|x|^p).$
The proof of the comparison result relies on classical techniques of
viscosity solutions. We build a suitable test-function and
prove some fine estimates on the various terms which appear, the
main difficulty consists in dealing with the unbounded control
terms.
\smallskip

In Section~\ref{sec:super}, we extend the comparison result in~\cite{dll06} 
for equations with $p>1$ growth conditions on the datas (instead of
quadratic growth) with the additional nonlinearity $f.$
 One motivation to add the nonlinearity $f$ 
comes from the BSDEs (where $s=\sigma$) since
the main application of Theorem~\ref{thm-unicbis} is the
uniqueness for the equation stated in~\cite{bh08} (see Example~\ref{expl-bdse}).
In this case, we consider Hamiltonians $H$ with
$\alpha$-bounded diffusion matrices $\sigma,$ so
we choose to replace~\eqref{hjb}
by the control independent PDE~\eqref{model} 
to simplify the exposition.
The control case does not present additional difficulties with respect 
to~\cite[Theorem 2.1]{dll06}.
The main difficulty in the proof of Theorem~\ref{thm-unicbis}
is to be able to deal with solutions
growing like $O(1+|x|^p)$
(which are not bounded neither from above nor from below).
The strategy of proof is similar to the one used in~\cite{dll06} 
which consists essentially in the following three steps. 
First one computes the equation satisfied by $w_{\mu}=\mu u-v$, 
being $u,v$ respectively the subsolution and the supersolution of the
original PDE and $0<\mu<1$ a parameter. Then for all $R>0$ one
constructs a strict supersolution $\Phi^R_\mu$ of the ``linearized
equation" such that $\Phi^R_\mu(x,t)\to 0$ as $R\to +\infty.$
Finally one shows that $w_{\mu}\le \Phi^R_\mu$ and one concludes 
by letting first $R\to +\infty$ and then $\mu\to 1.$
\smallskip

A by-product of the comparison results obtained in
Sections~\ref{sec:unbd} and~\ref{sec:super} and Perron's Method
  of Ishii~\cite{ishii87} is the existence and uniqueness of a
 continuous solution to~\eqref{hjb} which is respectively  
$o(1+|x|^p)$ and $O(1+|x|^p)$  as $|x|\to\infty.$ 
However, under our general assumptions one cannot expect
the existence of a solution for all times as Example
\ref{LP-deterministe} shows.
\smallskip

Let us compare our results with related ones in the literature
for such kind of Hamilton-Jacobi equations.
Uniqueness and existence problems for a class of first-order Hamiltonians corresponding to
unbounded control sets and under assumptions including deterministic linear qua\-dra\-tic
problems have been addressed by several authors, see, e.g. the book of Bensoussan~\cite{bensoussan82},
the papers of Alvarez~\cite{alvarez97}, Bardi and Da Lio~\cite{bdl97}, Cannarsa and Da Prato~\cite{cdp89},
Rampazzo and Sartori~\cite{rs00} in the case of convex operators, and the papers of
Da Lio and McEneaney~\cite{dlme02} and Ishii~\cite{ishii97} for more general operators.
As for second-order Hamiltonians under quadratic growth assumptions, Ito~\cite{ito01} obtained
the existence of locally Lipschitz solutions to particular equations of the form \eqref{hjbintr}
under more regularity conditions on the data, by establishing a priori estimates on the solutions.
Whereas Crandall and Lions in~\cite{cl90} proved a uniqueness result for very particular operators
depending only on the Hessian matrix of the solution. In the case of quasilinear degenerate parabolic equations,
existence and uniqueness results for viscosity solutions which may
have a quadratic growth are proved in~\cite{bbbl03}.
The results which are the closest to ours were obtained in the
following works.
Alvarez~\cite{alvarez96} addressed the case of stationary less
general equations (see Example~\ref{exple-simple-p}).
Krylov~\cite{krylov01} succeeded in dealing with equations
encompassing the classical Stochastic Linear Quadratic problem
but his assumptions are designed to handle exactly this case
(cf. Example~\ref{slqp} and the discussion therein).
Finally Kobylanski~\cite{kobylanski00} studied 
also~\eqref{hjb} under quite general assumptions on the datas
but for bounded solutions. It seems
to be difficult to obtain such a generality in the case of unbounded
solutions since her proof is based on changes of functions of the form
$u\to -e^{-u}$ which do not work for solutions which are neither
bounded from below nor from above.
\smallskip

The rigorous connection between control problems and
Hamilton-Jacobi-Bellman equations is not addressed in this paper.
In the  framework of unbounded controls it may be  rather delicate.
 Some results in this direction were obtained for infinite horizon in
 the deterministic case by Barles \cite{barles90b} and in the stochastic case by 
Alvarez~\cite{alvarez96, alvarez97}, Krylov~\cite{krylov01} and by the authors \cite{dll06}.
\smallskip

Finally, let us mention that the convexity of the operator with
respect to the gradient is crucial in our proofs.
The case of Hamiltonians which are  neither convex nor concave 
(which, in the case of Equations~\eqref{isaacs-eq}, amounts to take  
both the control sets $A$ and $B$ unbounded) is also of interest
and it is a widely open subject.
Some results in this direction were obtained 
in~\cite[Section 4]{dll06},
for instance in the case of first order equations of the form
 $$
\frac{\partial u}{\partial t} +h(x,t)|Du|^2=0~~~\mbox{in
  $\R^N\times[0,T],$}
$$
where $h(x,t)$ may change sign and $u$ has a quadratic growth.
In a forthcoming paper we are going to investigate this
issue for more general quadratic non convex-non concave equations.
\smallskip

Throughout the paper we will use the following notations. 
For all integer $N,M\geq 1$ we denote by ${\cal{M}}_{N, M} (\R)$  
(respectively ${\cal{S}}_{N} (\R),$ ${\cal{S}}_{N}^+ (\R)$)  the set of  real  $N\times M$ matrices  
(respectively real symmetric matrices, real symmetric nonnegative  $N\times N$ matrices).  
For the sake of notations, all the norms which appear in the sequel are denoting by $|\cdot|.$ 
The standard Euclidean inner product in $\R^N$ is written $\langle\cdot,\cdot\rangle.$  
We recall that a modulus of continuity $m: \R\to \R^+$ is a nondecreasing continuous function 
such that $m(0)=0.$ We set $B(0,R)= \{ x\in \R^N : |x|<R\}.$  Finally
for any $O\subseteq\R^K$, 
we denote by $USC(O)$ the set of 
upper semicontinuous functions in $O$ and by $LSC(O)$ the set of 
lower semicontinuous functions in $O.$    Given $p>1$ we will denote by $p^\prime$ its coniugate, namely
$$
\frac{1}{p}+ \frac{1}{p'}=1.
$$

\noindent{\bf Acknowledgments.} Part of this work was done while the
second author was a visitor at the FIM at the ETH in Z\"urich in
January 2007. He would like
to thank the Department of Mathematics for his support.
We thank Guy Barles for useful comments on the first version
of this paper.

\section{Hamilton-Jacobi-Bellman equations with unbounded diffusion in the control  } 
\label{sec:unbd} 

In this  Section we  prove a comparison result for second-order 
fully nonlinear partial differential equations of the form~\eqref{hjb}.
 The main difference with respect to the result in~\cite{dll06} 
is that here we suppose that the diffusion matrix $\sigma$ 
depends in a unbounded way in the control (see condition~\eqref{sig-non-borne}).
  The compensation to the condition~\eqref{sig-non-borne} is that we are able to get the uniqueness
 result in the smaller  class of functions which are $o(1+|x|^p)$ as
 $|x|\to\infty$ (see~\eqref{classeC}).  

  We list below the main  assumptions on $H$ and $f$.   \par
\par\medskip 
{\noindent \textbf{(A)} (Assumption on $H$)~:}\par 
\begin{enumerate} 
\item[(i)]~$A$ is a subset of a separable complete normed space. 
The main point here is the possible unboundedness of $A.$ 

\item[(ii)]~${b}\in C(\R^N \times [0,T]\times A ; \R^N)$ and
    there exists $C_b>0$ such that, for all 
$x,y\in\R^N,$ $t\in [0,T],$ ${\alpha}\in A,$ 
\begin{eqnarray*} 
|{b}(x,t,{\alpha})-{b}(y,t,{\alpha})|& \leq &{C_b} (1+|{\alpha}|)|x-y|, \\ 
|{b}(x,t,{\alpha})| &\leq & {C_b} (1+|x|+|{\alpha}|)\ ; 
\end{eqnarray*} 
\item[(iii)]~ 
${\ell}\in C(\R^N \times [0,T]\times A ; \R)$ and, there exist
  $p>1$ and $ C_\ell, \nu>0$ such that, for all $x\in\R^N,$ $t\in
  [0,T],$ ${\alpha}\in A,$ 
$$ 
C_\ell (1+|x|^p+|{\alpha}|^p) 
\geq 
{\ell}(x,t,{\alpha})
\geq {\nu} |{\alpha}|^p  -C_{\ell} (1+|x|^p)
$$ 
and for every $R>0,$ there exists a modulus of continuity $m_R$ 
such that  for all $  x,y\in B(0,R),$ $t\in [0,T],$ ${\alpha}\in A,$ 
\begin{equation*} 
 |{\ell}(x,t,{\alpha})-{\ell}(y,t,{\alpha})| \leq (1+|{\alpha}|^p)\, m_R(|x-y|)\ ; 
\end{equation*} 
\item[(iv)] 
${\sigma}\in C(\R^N\times [0,T]\times A ; {\cal{M}}_{N,M} (\R))$ is 
Lipschitz continuous with respect to $x$ with a constant independent of $(t, 
{\alpha})$: namely, there exists $C_\sigma >0$ such that, 
for all $x,y\in \R^N$ and $(t, {\alpha}) \in [0,T]\times A,$ 
$$ 
|{\sigma}(x,t, {\alpha})- {\sigma}(y,t, {\alpha})|\leq C_\sigma |x-y|, 
$$ 
and satisfies { for \ every} $  x\in \R^N, \ t\in [0,T], \ {\alpha}\in
A,$ 
\begin{eqnarray}  \label{sig-non-borne}
|{\sigma}(x,t, {\alpha})| \leq  {C_{\sigma}}(1+|x|+|\alpha|). 
\end{eqnarray}
\end{enumerate}

\noindent{\bf (B)}  (Assumption on $f$)\\
$f\in C([0,T]\times \R^N\times\R \times \R^N; \R)$
and, for all $R>0,$ there
exist a modulus of continuity $m_R$ and $C_s, \hat{C}>0$ such that,
for all $t\in [0,T],$ $x,y\in\R^N,$ $u,v\in\R,$ $z\in\R^N,$ 
\begin{eqnarray*} 
&{\rm (i)}& |f(x,t,u,z)|\leq C_f (1+|x|^{p } +|u|+|z|^{p^\prime}), \\
&{\rm (ii)}& |f(x,t,u,z)-f(y,t,u,z)|\leq  m_R( (1+|u|+|z|)|x-y|) \ \ \ {\rm
    if} \ |x|+|y|\leq R,\\
&{\rm (iii)}& z\mapsto f(x,t,u,z) \ {\rm is \ convex},\\
&{\rm (iv)}& s\in C(\R^N\times [0,T]; \mathcal{M}_N), \ \ \
          |s(x,t)-s(y,t)|\leq C_s |x-y|, \ \ \ |s(x,t)|\leq C_s,\\
&{\rm (v)}&|f(x,t,u,z)-f(x,t,v,z)|\leq \hat{C} |u-v|. 
\end{eqnarray*} 

The typical case we have in mind in the context of {\textbf{(A)}}(iv)
($\sigma$ not bounded with respect to the control) is
\begin{eqnarray*} 
\sigma (x,t,\alpha)= Q(t)x + R(t)\alpha,
\end{eqnarray*} 
where $Q(t)$ and $R(t)$ are matrices of suitable sizes. This case
includes Linear Quadratic control problems, see Example~\ref{slqp}.

Under the current hypotheses, the Hamiltonian $H$ may be 
infinite (see Example~\ref{slqp})
and for this reason we re-formulate the definition of viscosity solution in the following way.

\begin{Definition}\label{newdef} \ \\ 
(i) A function $u\in USC(\R^N\times [0,T])$  is a viscosity 
subsolution of~\eqref{hjb} if for all $(x,t)\in\R^N\times[0,T]$ and $\varphi 
\in C^2(\R^N\times [0,T])$ such that $u-\varphi$ has a maximum at $(x,t),$ 
we have $u(x,t)\leq \psi(x)$ if $t=0$ and, if $t>0,$ then 
$$ 
\frac{\partial \varphi}{\partial t}(x,t)+H(x,t,D\varphi (x,t), D^2\varphi 
(x,t))+f(x,t,u(x,t),s(x,t)D\varphi (x,t))\leq 0, 
$$ 
which is equivalent to: for all $\alpha \in A,$ 
\begin{eqnarray} \label{def-sub}
&&\frac{\partial \varphi}{\partial t}(x,t)- 
\langle {b}(x,t,{\alpha}),  D\varphi(x,t)\rangle - {\ell}(x,t,{\alpha}) 
-{\rm Trace}\left[ {\sigma}(x,t, {\alpha}){\sigma}^T(x,t,{\alpha}) 
 D^2\varphi(x,t) \right]\nonumber\\
 && ~~~~+f(x,t,u(x,t),s(x,t)D\varphi (x,t))\leq 0. 
\end{eqnarray}
(ii)
A function $u\in USC(\R^N\times [0,T])$  is a viscosity 
supersolution of~\eqref{hjb} if for all $(x,t)\in\R^N\times[0,T]$ and $\varphi 
\in C^2(\R^N\times [0,T])$ such that $u-\varphi$ has a minimum at $(x,t),$ 
we have $u(x,t)\geq \psi(x)$ if $t=0$ and, if $t>0,$ then 
 for all $\eta >0,$ there exists $\alpha_\eta =\alpha(\eta,x,t)\in A,$ such that 
\begin{eqnarray} \label{def-super}
&&\frac{\partial \varphi}{\partial t}(x,t)- 
\langle {b}(x,t,{\alpha_\eta}),  D\varphi(x,t)\rangle - {\ell}(x,t,{\alpha_\eta}) 
-{\rm Trace}\left[ {\sigma}(x,t, {\alpha_\eta}){\sigma}^T(x,t,{\alpha_\eta}) 
 D^2\varphi(x,t) \right]\nonumber\\
 &&~~~~ +f(x,t,u(x,t),s(x,t)D\varphi (x,t))\geq -\eta. 
\end{eqnarray}
(iii) A locally bounded function $u: \R^N\times [0,T]\to \R$ is a viscosity 
solution 
of~\eqref{hjb} if its USC envelope $u^*$ is a subsolution and its LSC 
envelope $u_*$ 
is a supersolution. 
 \end{Definition} 
Note that~\eqref{def-sub} and~\eqref{def-super} is only a way to write the definition
of sub- and supersolutions without writing a supremum which could not exist because of
assumption~\eqref{sig-non-borne}.


We say that a function $u:\R^N\times [0,T]\to \R$ is in the class 
$\mathcal{C}_p$ 
if 
\begin{eqnarray} \label{classeC}
\frac{u(x,t)}{1+|x|^p}\ \mathop{\longrightarrow}_{|x|\to +\infty} \ 0, 
\quad \text{uniformly with respect to $t\in  [0,T].$} 
\end{eqnarray} 
 Note that $u\in \mathcal{C}_p$ if and only if, for all
$\eps>0,$ there exists $M_\eps >0$ such that
\begin{eqnarray*} 
|u(x,t)|\leq M_\eps + \eps (1+|x|^p) \quad \text{for all $(x,t)\in
  \R^N\times [0,T].$}
\end{eqnarray*} 
 In particular, for all $\lambda >0,$ 
\begin{eqnarray} \label{maxC} 
\mathop{\rm sup}_{x\in \R^N} \{ u(x,t)-\lambda (1+|x|^p)\}= M_\lambda < +\infty.
\end{eqnarray} 

 \par The main result of this Section is 
the 

 \begin{Theorem} \label{thm-unicite} 
Assume {\textbf{(A)}}-{\textbf{(B)}} and suppose that $\psi$ is a continuous function which
belongs to $\mathcal{C}_p.$ 
Let $u \in USC(\R^N\times [0,T])$ be a viscosity subsolution of~\eqref{hjb} 
and $v \in LSC(\R^N\times [0,T])$ be a viscosity supersolution of~\eqref{hjb}. 
Suppose that $U$ and $V$ are in the class $\mathcal{C}_p$ defined 
by~\eqref{classeC} and  
satisfy $u(x,0)\leq\psi (x)\leq v(x,0).$  
Then $u\leq v$ in $\R^N\times [0,T].$ 
\end{Theorem} 

Before giving the proof of the theorem, let us state an existence
result and some examples of applications.
As it was already observed in~\cite{dll06}, the question of the
existence of a continuous solution to~\eqref{hjbintr} is
not completely obvious and in general the solutions may exist
only for short time (see Example~\ref{LP-deterministe}). One way to obtain the
existence is to establish a link between the solution of the PDE
and related control problems or BDSE systems which have a solution.
By using PDE methods,
in the framework of viscosity solutions, the existence
is usually  a consequence of the comparison principle by
means of Perron's method, as soon as we can build a sub- and a super-solution
to the problem. Here, the comparison principle is proved in the class of functions
belonging to $\mathcal{C}_p.$ 
Therefore, to prove the existence, it siffices to
build   sub- and super-solutions to~\eqref{model} in   $\mathcal{C}_p.$
We need to strengthen {\bf (A)}(iii)
and {\bf (B)}(i) by assuming that $\ell(\cdot, t,\alpha), f(\cdot,t,u,z)\in\mathcal{C}_p$
uniformly with respect to $\alpha, t,u,z,$ i.e., for all
$(x,t,\alpha,u,z)\in \R^N\times [0,T]\times A \times \R\times \R^N,$
\begin{eqnarray}\label{assumpt-sup} 
\begin{array}{c}
\chi(x)\geq  \ell (x,t,\alpha)\geq \nu |\alpha|^p -\chi(x),
\quad |f(x,t,u,z)|\leq C_f (1+\gamma(x) +|u|+|z|^{p^\prime}),\\
\displaystyle
{\rm and}\quad \mathop{\rm lim}_{|x|\to +\infty} \frac{\chi (x)}{1+|x|^p},\,
\frac{\gamma(x)}{1+|x|^p} =0.
\end{array}
 \end{eqnarray} 
We have

\begin{Theorem} \label{thm-existence2}
Assume \textbf{(A)}--\textbf{(B)} and~\eqref{assumpt-sup}. 
For all $\psi\in \mathcal{C}_p,$  there is $\tau >0$ such that
there exist a subsolution $\underline u\in \mathcal{C}_p$
and a supersolution $\overline u\in \mathcal{C}_p$
of~\eqref{hjb} in $\R^N\times [0,T].$
In consequence, Equation~\eqref{hjb} 
 has a unique continuous viscosity solution
in $\R^N\times [0,\tau ]$  in the class $\mathcal{C}_p.$
\end{Theorem}

The proof of this theorem is postponed at the end of the section.

\begin{Example} [A Stochastic Linear Quadratic Control Problem] \label{slqp}
Consider the stochastic differential equation
(in dimension 1 for sake of simplicity)
$$
\left\{\begin{array}{ll}
dX_s=X_sds + \sqrt{2}\alpha_s dW_s,&  t\leq s\leq T, \ t\in(0,T],\\
X_0=x\in\R,
\end{array}\right.
$$
where $W_s$ is a standard Brownian motion, $(\alpha_s)_s$ is a real valued progressively measurable process and
the value function is given by
$$
V(x,t)=\mathop{\rm inf}_{(\alpha_s)_s}E_{tx}\left\{ \int_t^T|\alpha_s|^2\,ds+\psi(X_T)\right\}.
$$
(Note that in this case, $p=p'=2$.)
The Hamilton-Jacobi equation formally associated to this problem is
\begin{eqnarray}\label{hj-lq} 
\left\{
\begin{array}{ll}
\displaystyle
-u_t+\sup_{\alpha\in \R}\{-\alpha^2(u''+1)\}- x u' =0  & \mbox{in $\R\times(0,T],$}\\
u(x,T)=\psi(x) & \mbox{in $\R$}.
\end{array}\right.
\end{eqnarray}
We observe    that in this case if $u''+1<0$ then the Hamiltonian becomes $+\infty$.
Nevertheless, we are able to prove comparison~\eqref{hj-lq} as soon as
 the terminal cost $\psi\in \mathcal{C}_p$ (i.e., has a strictly sub-$p$ growth).
This is not completely satisfactory since, in the classical
Linear Quadratic Control Problem, one expects to have quadratic
terminal costs like $\psi(x)=|x|^2.$
Let us mention that Krylov~\cite{krylov01} succeeded in treating
this latter case. But his proof consists on some algebraic
computations which rely heavily on the particular form of the datas
(the datas are supposed to be polynomials of degree 1 or degree 2 in $(x,\alpha)$).
In our case, up to restrict slightly the growth, we are able to deal
with general datas.
\end{Example}  

\begin{Remark} \label{ajout-G}
Theorems~\ref{thm-unicite} and~\ref{thm-unicbis} still hold for the
Isaacs equation of~\cite{dll06},
\begin{eqnarray} \label{isaacs-eq}
\displaystyle\frac{\partial u}{\partial t}
+ H(x,t,Du, D^2 u)+ G(x,t,Du, D^2 u) + f(x,t,u,s(x,t)Du) = 0 
\end{eqnarray}
where
\begin{eqnarray*}
G(x,t, q,X) =
\mathop{\rm inf}_{\beta\in B}
\left\{ -\langle g(x,t, \beta ),  q \rangle
- l (x,t, \beta) -{\rm Trace}\left[ c(x,t, \beta)c^T(x,t,\beta)X \right] \right\},
\end{eqnarray*}
is a concave Hamiltonian, 
$B$ is bounded, $g,l,c$ satisfy respectively {\bf (A)}(ii),(iii),(iv)
(with bounded controls $\beta$).
The case where both the control sets $A$
and $B$ are unbounded is rather delicate. It is the aim of a future work.
\end{Remark}  

Let us turn to the proof of the comparison theorem.\\

\noindent{\bf Proof of Theorem \ref{thm-unicite}.} We are going to show that for every  
  $\mu\in(0,1),$ $\mu u-v\le 0,$ in $\R^N\times[0,T].$  To this end 
we argue by contradiction assuming that there exists
$(\hat{x},\hat{t})\in\R^N\times [0,T]$ such that
\begin{eqnarray}\label{contr879} 
 u(\hat{x},\hat{t})- v(\hat{y},\hat{t})>\delta>0.
\end{eqnarray} 
We divide the proof  in several steps. \\

\noindent{\it 1. The $\mu$-equation for the subsolution.}
If $u$ is a subsolution of~\eqref{hjb}, then $\tilde{u}=\mu u$
is a subsolution of
\begin{eqnarray*}
&&\tilde{u}_t+\sup_{\alpha\in A}\{- {\rm Trace}\big(\sigma(x,t,\alpha)\sigma(x,t,\alpha)^T D^2 \tilde{u}\big) 
+\langle b(x,t,\alpha),D\tilde{u}\rangle-\mu\ell(x,t,\alpha)\}\\
&& +\mu f\big(x,t,\frac{1}{\mu} \tilde{u}(x,t),\frac{1}{\mu}s(x,t)D\tilde{u})\le 0, 
\end{eqnarray*}
with the initial condition $\mu u(x,0)\leq \mu\psi(x).$ \\

\noindent{\it 2. Test-function and estimates on the penalization terms.}
For all $\eps>0$, $\eta>0$ and $\theta, L>0$ (to be chosen later) we consider 
  the auxiliary function 
  \begin{eqnarray*}\label{auxfunctbis} 
  \Phi(x,y,t)=\mu u(x,t)-v(y,t)-e^{Lt}\big( \frac{|x-y|^2}{\eps^2}+ \theta(1-\mu)(1+|x|^2+|y|^2)^{p/2}\big)
-\rho t. 
  \end{eqnarray*} 
Since $u,v\in \mathcal{C}_{p},$ the supremum of  $\Phi$ in 
$\R^N\times\R^N\times[0,T]$ is achieved at a point  $(\bar x,\bar y,\bar t).$ 
We will drop for simplicity of notation the dependence 
  on the various parameters. If $\theta$ and $\rho$ are small enough we have
  \begin{eqnarray*} 
   \Phi(\bar x,\bar y,\bar t)\ge \mu u(\hat x,\hat y)-v(\hat x,\hat y)  
-\theta(1-\mu)(1+2|\hat x|^p)-\rho \hat t>\frac{\delta}{2},
   \end{eqnarray*} 
which implies
\begin{eqnarray*} 
 \frac{|\bar x-\bar y|^2}{\eps^2}+ \theta(1-\mu)(1+|\bar x|^2+|\bar y|^2)^{p/2}
\leq
\mu u(\bar x,\bar t)-v(\bar y,\bar t).
\end{eqnarray*} 
Therefore, by~\eqref{maxC}, we get 
\begin{eqnarray*} 
&& \frac{|\bar x-\bar y|^2}{\eps^2}+ \theta\frac{1-\mu}{2}(1+|\bar x|^2+|\bar y|^2)^{p/2}\\
&\leq & 
\mathop{\rm sup}_{(x,t)\in\R^N\times [0,T]}\{ \mu u(x,t)-\theta\frac{1-\mu}{2}(1+| x|^{p})\}\\
&&
+ \mathop{\rm sup}_{(x,t)\in\R^N\times [0,T]} \{ -v(x, t)-\theta\frac{1-\mu}{2}(1+| x|^{p})\}\\
&\leq & 
M 
\end{eqnarray*} 
for some $0<M=M(\mu,\theta,u,v).$
Thus 
\begin{eqnarray} \label{xybornes}
|\bar x|, |\bar y|\leq R_{\mu,\theta}
\end{eqnarray}
 with $R_{\mu,\theta}$ independent of
$\varepsilon$ and $|\bar x-\bar y|\to 0$ as $\varepsilon\to 0.$ 
Up to extract a subsequence, we can assume that
\begin{eqnarray} \label{def-compxy}
\bar x, \bar y \to x_0\in \overline{B}(0,R_{\mu,\theta}),~~\bar t\to t_0~~\mbox{as $\eps\to 0$}
\end{eqnarray}
Actually we can obtain a more precise estimate: we have
 \begin{eqnarray*} 
 \Phi(\bar x,\bar y,\bar t)\geq \mathop{\rm max}_{\R^N\times [0,T]}
\{  \mu u(x,t) -v(x,t)- e^{Lt}\theta(1-\mu)(1+2|x|^2)^{p/2}-\rho t \}:=  M_{\mu,  \theta}.
 \end{eqnarray*} 
Thus
 \begin{eqnarray*} 
\mathop{\rm lim\,inf}_{\eps\to 0}\Phi(\bar x,\bar y,\bar t)\geq M_{\mu,  \theta}.
 \end{eqnarray*} 
On the other hand
 \begin{eqnarray*} 
   & & 
   \limsup_{\eps\to 0}\Phi(\bar x,\bar y,\bar t) \\ 
& & \le  
\limsup_{ \eps\to 0}\, [ \mu u(\bar x,\bar t)\!-\!v(\bar y,\bar t) 
-e^{L\bar t}\theta(1\!-\!\mu)(1+|\bar x|^2+|\bar y|^2)^{p/2}\!-\!\rho \bar t] 
  -\liminf_{\eps\to 0}e^{L\bar t}\frac{|\bar x\!-\!\bar y|^2}{\eps^2}  \\ 
& & \le  M_{\mu,  \theta}-\liminf_{\eps\to 0}e^{L\bar t}\frac{|\bar x-\bar y|^2}{\eps^2}. 
   \end{eqnarray*} 
   By combining the above inequalities we get, up to subsequences, that 
\begin{eqnarray} 
\label{def-compxy2}
  && \frac{|\bar x-\bar y|^2}{\eps^2}\to 0 ~~ \mbox{as $\eps\to 0.$} 
 \end{eqnarray} 
Note that we have
\begin{eqnarray} \label{mod45}
|\bar x-\bar y|, \, \frac{|\bar x-\bar y|^2}{\eps^2}= m(\eps),
 \end{eqnarray} 
where $m$ denotes a modulus of continuity independent of $\varepsilon$
(but which depends on $\theta, \mu$). \\

\noindent{\it 3. Ishii matricial theorem and viscosity inequalities.}
We set 
  $$ 
  \Theta(x,y,t)=e^{Lt}\big( \frac{|x-y|^2}{\eps^2}+ \theta(1-\mu)(1+|x|^2+|y|^2)^{p/2}\big)+\rho t. 
  $$ 
We claim that there is a subsequence $\eps_n$ such that $\bar t=0.$
Suppose by contradiction that for all $\eps>0$ we have $\bar t>0$.
Next Steps are devoted to prove some estimates in order to obtain
the desired contradiction at the end of Step 8.

By  Theorem 8.3 in the User's guide~\cite{cil92}, for 
every $\varrho >0,$ there exist 
$a_1, a_2 \in\R$ and $X,Y\in{\mathcal{S}}_{N}$ such that 
\begin{eqnarray*} 
\left(a_1,D_x\Theta(\bar x,\bar y,\bar t), 
X\right) 
  \in\bar{\mathcal{P}}^{2,+}(\mu u)(\bar x, \bar t), \\[2mm] 
\left(a_2,-D_y \Theta(\bar x,\bar y,\bar t), 
Y\right) 
  \in\bar{\mathcal{P}}^{2,-}(v)(\bar y,\bar t), \\[2mm] 
a_1-a_2=\Theta_t (\bar x,\bar y,\bar t),
\end{eqnarray*} 
and 
\begin{eqnarray*}\label{Ineq-Matricebis} 
-(\frac{1}{\varrho}+|{{M}}|)I \leq 
\left( 
\begin{array}{cc} 
X & 0\\ 
0 & -Y\\ 
\end{array} 
\right)\leq {{M}}+ \varrho {{M}}^2 
 \end{eqnarray*} 
where 
${{M}}=D^2\Theta (\bar x,\bar y,\bar t ).$ 
Note that 
$$ 
a_1-a_2=Le^{L\bar t} \big( \frac{|\bar x-\bar y|^2}{\eps^2}
+ \theta(1-\mu)(1+|\bar x|^2+|\bar y|^2)^{p/2}\big) +\rho, 
$$ 
and, setting $\displaystyle{p_\varepsilon = 2 e^{L\bar t} \frac{\bar x-\bar y}{\varepsilon^2}},$ 
$q_x=e^{L\bar t} p \theta(1-\mu)\bar x(1+|\bar x|^2+|\bar y|^2)^{p/2-1}$,  
$q_y=-e^{L\bar t}p \theta(1-\mu)\bar y (1+|\bar x|^2+|\bar y|^2)^{p/2-1}$ we have 
$$ 
D_x\Theta(\bar x,\bar y,\bar t) 
=  p_\varepsilon + q_x 
\ \ \ {\rm and} \ \ \ 
 D_y\Theta(\bar x,\bar y,\bar t) 
= -p_\varepsilon-q_y, 
$$ 
and 
$$ 
M=A_1+A_2+A_3$$ 
where 

$$
A_1= \frac{2e^{L\bar t}}{\varepsilon^2}\left( 
\begin{array}{cc} 
 I & -I\\ 
-I & I\\ 
\end{array} 
\right),
$$ 
$$
A_2= e^{L\bar t}p \theta(1-\mu)(1+|\bar x|^2+|\bar y|^2)^{p/2-1}\left( 
\begin{array}{cc}   
 I 
& 0\\ 
0 & I \\ 
\end{array} 
\right),
$$ 
$$
A_3= e^{L\bar t}p (p-2)\theta(1-\mu)(1+|\bar x|^2+|\bar y|^2)^{p/2-2}\left( 
\begin{array}{cc}   
 x\otimes x 
& x\otimes y\\ 
x\otimes y & y\otimes y \\ 
\end{array} 
\right). 
$$ 
It follows 
\begin{eqnarray}
\langle X \xi,\xi\rangle - \langle Y \zeta,\zeta\rangle &\leq  &  
 \frac{2e^{L\bar t}}{\varepsilon^2}   |\xi-\zeta|^2 \nonumber\\ &&+e^{L\bar t}p \theta(1-\mu)(1+|\bar x|^2+|\bar y|^2)^{p/2-1}(|\xi|^2+|\zeta|^2) \nonumber\\ 
&&+ 2 e^{L\bar t}p (p-2)\theta(1-\mu)(1+|\bar x|^2+|\bar
y|^2)^{p/2-2}\left(\langle\xi,x\rangle^2+\langle\zeta,y\rangle^2\right)\nonumber\\
&&+ 
m\left(\frac{\varrho}{\varepsilon^4}\right),\label{ineg-mat555}
\end{eqnarray}
where $m$ is a modulus of continuity which is independent of 
$\rho$ and $\varepsilon.$ 

We now write the viscosity inequalities satisfied by the subsolution $\mu u$
and the supersolution $v$  (recall that we assume $\bar t>0$).

For all $\alpha\in A$ we have  
\begin{eqnarray}\label{subineq}
 && a_1-\mbox{Trace}(\sigma(\bar x,\bar t,\alpha)\sigma(\bar x,\bar t,\alpha)^T X)
+\langle b(\bar x,\bar t,\alpha), p_\varepsilon+q_x \rangle
  -\mu\ell(\bar x,\bar t,\alpha)\nonumber\\&&~~+\mu f(\bar x,\bar t,u(\bar x,\bar t ),\frac{1}{\mu}s(\bar x,\bar t)(p_\eps+q_x))\le 0. 
\end{eqnarray} 
On the other hand, for all $\eta>0,$ there exists $\alpha_\eta\in A$ such that 
\begin{eqnarray}\label{superineq}
 && a_2-\mbox{Trace}(\sigma (\bar y,\bar t,\alpha_{\eta})\sigma^T(\bar y,\bar t,\alpha_{\eta}) Y)
+\langle b(\bar y,\bar t,\alpha_{\eta}), p_\varepsilon+q_y \rangle
-\ell(\bar y,\bar t,\alpha_{\eta})\nonumber\\&&~~+f(\bar y,\bar t,v(\bar y,\bar t ),\frac{1}{\mu}s(\bar y,\bar t)(p_\eps+q_y))\geq -\eta.
\end{eqnarray}
We set for simplicity
\begin{eqnarray*}
&&\sigma_x:={\sigma}(\bar x,\bar t,\alpha_{\eta}),\
\sigma_y={\sigma}(\bar y,\bar t,\alpha_{\eta}) \\
&& b_x= {b}(\bar x,\bar t,\alpha_{\eta}), \
b_y= {b}(\bar y,\bar t,\alpha_{\eta}),\ s_x=s(\bar x,\bar t),\ s_y=s(\bar y,\bar t).
 \end{eqnarray*}
By subtracting~\eqref{subineq} and~\eqref{superineq} we get
\begin{eqnarray}\label{ineg-651} 
&& Le^{L\bar t} \big( \frac{|\bar x-\bar y|^2}{\eps^2}
+ \theta(1-\mu)(1+|\bar x|^2+|\bar y|^2)^{p/2}\big) +\rho\nonumber\\
&\leq & 
\mbox{Trace}(\sigma_x \sigma_x^T X- \sigma_y \sigma_y^T Y)
+ \langle b_y, p_\varepsilon+q_y \rangle- \langle b_x, p_\varepsilon+q_x \rangle\nonumber\\
&&  -\ell(\bar y,\bar t,\alpha_{\eta})+\mu \ell(\bar x,\bar t,\alpha_{\eta})\nonumber\\
&& + f(\bar y,\bar t,v(\bar y,\bar t ),s_y(p_\eps+q_y))-\mu f(\bar x,\bar t,u(\bar x,\bar t ),\frac{1}{\mu}s_x(p_\eps+q_x))+\eta.
\end{eqnarray}

 \noindent{\it 4. Estimates of the second-order terms.}
From~\eqref{ineg-mat555} and {\bf (A)}(iv), it follows 
\begin{eqnarray*}
 & & {\rm Trace}\left[ {\sigma_x} {\sigma_x}^T X -{\sigma_y} {\sigma_y}^T Y \right] 
- m\left(\frac{\varrho}{\varepsilon^4}\right)
\nonumber\\
&\leq & 
e^{L\bar t}\left(
\frac{2}{\eps^2}|\sigma_x-\sigma_y|^2
+ p \theta(1-\mu)(1+|\bar x|^2+|\bar y|^2)^{p/2-1}(|\sigma_x|^2+|\sigma_y|^2)\right.\nonumber\\
&&\hspace*{1.2cm}\left. +2 p (p-2)\theta(1-\mu)(1+|\bar x|^2+|\bar y|^2)^{p/2-2}
(|\sigma_x|^2|x|^2+|\sigma_y|^2|y|^2)\rule{0cm}{0.6cm}\right)
\nonumber\\
& \leq & 2C_\sigma^2e^{L\bar t}\left(
\frac{|\bar x-\bar y|^2}{\eps^2} +
  p(p-1) \theta(1-\mu)(1+|\bar x|^2+|\bar y|^2)^{p/2-1}(1+|\bar x|^2+|\bar y|^2+|\alpha_\eta|^2)\right)\nonumber\\
&\leq &
2C_\sigma^2e^{L\bar t}\left(
\frac{|\bar x-\bar y|^2}{\eps^2} +
p(p-1) \theta(1-\mu)(1+|\bar x|^2+|\bar y|^2)^{p/2}\right.\nonumber\\
&&\hspace*{2.4cm}\left.
+p(p-1) \theta(1-\mu) |\alpha_\eta|^2(1+|\bar x|^2+|\bar y|^2)^{p/2-1}
\rule{0cm}{0.6cm}\right).
\nonumber 
\end{eqnarray*} 
By Young's inequality,
\begin{eqnarray*}
|\alpha_\eta|^2(1+|\bar x|^2+|\bar y|^2)^{p/2-1}
\leq \frac{2}{p}|\alpha_\eta|^{p}+\frac{p-2}{p}(1+|\bar x|^2+|\bar y|^2)^{p/2}.
\end{eqnarray*}
It follows, using~\eqref{mod45},
\begin{eqnarray}
{\rm Trace}\left[ {\sigma_x} {\sigma_x}^T X -{\sigma_y} {\sigma_y}^T Y \right] 
&\leq&
4(p-1)^2C_\sigma^2e^{L\bar t}\theta(1-\mu)(1+|\bar x|^2+|\bar y|^2)^{p/2}\nonumber\\
&&\hspace*{0.6cm}
+ 4(p-1)C_\sigma^2e^{L\bar t}\theta(1-\mu)|\alpha_\eta|^{p}
+ m(\eps)+m\left(\frac{\varrho}{\varepsilon^4}\right)\!. \label{ineg222}
\end{eqnarray}

\noindent{\it 5. Estimates of the drift terms.}
By  using  {\bf (A)}(ii) and, from \eqref{mod45}, 
by taking  $\eps$ is small enough in order that $|\bar x-\bar y|\leq 1,$ we get 
\begin{eqnarray*}
&& \langle b_y, p_\varepsilon+q_y \rangle- \langle b_x, p_\varepsilon+q_x \rangle\nonumber\\
&\leq & 
 \langle b_y-b_x , p_\varepsilon+q_y \rangle +\langle b_x, q_y-q_x \rangle\nonumber\\
& \leq &
|b_y-b_x||p_\varepsilon|+ |b_y-b_x||q_y| + |b_x||q_x-q_y|\nonumber\\
&\leq & C_be^{L\bar t}\left(2(1+|\alpha_\eta|)\frac{|\bar x-\bar y|^2}{\eps^2}
+ 2{p}\theta(1-\mu)(1+|\alpha_\eta|)|\bar x-\bar y|(1+|\bar x|^2+|\bar y|^2)^{(p-1)/2}\right. \nonumber\\
&&\hspace*{1.5cm} \left. +2{p}\theta(1-\mu)(1+|\bar x|^2+|\bar y|^2)^{p/2}\rule{0cm}{0.6cm}\right)\nonumber\\
&\leq & C_be^{L\bar t}\left( 2\frac{|\bar x-\bar y|^2}{\eps^2}
+ 4{p}\theta(1-\mu)(1+|\bar x|^2+|\bar y|^2)^{p/2}\right. \nonumber\\
&&\hspace*{1.5cm} \left. +m(\eps)|\alpha_\eta|
+  \theta(1-\mu)m(\eps)|\alpha_\eta|(1+|\bar x|^2+|\bar y|^2)^{(p-1)/2}\rule{0cm}{0.6cm}\right). \nonumber
\end{eqnarray*} 
By Young's inequality, we get
\begin{eqnarray*}
&& m(\eps)|\alpha_\eta|+ \theta(1-\mu)m(\eps)|\alpha_\eta|(1+|\bar x|^2+|\bar y|^2)^{(p-1)/2}\\
&\leq&
\frac{ m(\eps)}{(\theta(1-\mu))^{1/(p-1)}}+\theta(1-\mu)|\alpha_\eta|^p
+ \theta(1-\mu)m(\eps) + \theta(1-\mu)(1+|\bar x|^2+|\bar y|^2)^{p/2}.
\end{eqnarray*} 
It follows
\begin{eqnarray}
&& \langle b_y, p_\varepsilon+q_y \rangle- \langle b_x, p_\varepsilon+q_x \rangle\nonumber\\
&\leq &
(4{p}+1)C_be^{L\bar t}\theta(1-\mu)(1+|\bar x|^2+|\bar y|^2)^{p/2}
+ C_be^{L\bar t}\theta(1-\mu)|\alpha_\eta|^p+ m(\eps).
\label{ineg-b}
\end{eqnarray} 

\noindent{\it 6. Estimates of running cost terms.}
Recall that we chose $\eps$ small enough in order that $|\bar x -\bar y|\leq 1.$
Setting $R=1,$ from {\bf (A)}(iii), we get 
 \begin{eqnarray*} 
 \mu {\ell}(\bar x,\bar t,{\alpha_\eta})- {\ell}(\bar y,\bar t,{\alpha_\eta})
&=&  (\mu -1)  {\ell}(\bar x,\bar t,{\alpha_\eta}) 
+  {\ell}(\bar x,\bar t,{\alpha_\eta})- {\ell}(\bar y,\bar t,{\alpha_\eta}) \nonumber\\ 
&\leq & 
(1-\mu) |\alpha_\eta|^{p} \left(-\nu + \frac{m_1(|\bar x-\bar y|)}{1-\mu}\right) \\ 
&& +  
C_\ell (1-\mu)(1+|\bar x|^p) +  m_1(|\bar x-\bar y|). \nonumber
 \end{eqnarray*}
Since $m_1(|\bar x-\bar y|)=m(\eps)$ by~\eqref{mod45}, we obtain
\begin{eqnarray} \label{ineg-ell}
  \mu {\ell}(\bar x,\bar t,{\alpha_\eta})- {\ell}(\bar y,\bar t,{\alpha_\eta})
\leq
(1-\mu) |\alpha_\eta|^{p} \left(-\nu + m(\eps)\right)
+ C_\ell (1-\mu)(1+|\bar x|^p) +  m(\eps).
 \end{eqnarray}
Note that it is the term ``$-(1-\mu) \nu|\alpha_\eta|^{p}$''
which will allow to control all the unbounded control terms in the sequel.\\

\noindent{\it 7. Estimates of $f$-terms.} To simplify, we replace
{\bf (B)}(v) by the assumption
that $f$ is nondecreasing with respect to the $u$ variable. By some
changes of functions as in Lemma~\ref{chgt-fct}, we can reduce to this
case without loss of generality.

We write
\begin{eqnarray*}
 f\big(\bar y,\bar t, v(\bar y,\bar t), s_y(p_\varepsilon+q_y)\big)
  - \mu\, f\big(\bar x,\bar t, u(\bar x,\bar t),  \frac{1}{\mu}s_x(p_\varepsilon+q_x)\big)
= \mathcal{T}_1 + \mathcal{T}_2 + \mathcal{T}_3 
\end{eqnarray*} 
with
\begin{eqnarray*}
&&  \mathcal{T}_1=
f\big(\bar y,\bar t, v(\bar y,\bar t), s_y(p_\varepsilon+q_y)\big)
- f\big(\bar x,\bar t, v(\bar y,\bar t), s_y(p_\varepsilon+q_y)\big),\\
&&  \mathcal{T}_2=
 f\big(\bar x,\bar t, v(\bar y,\bar t), s_y(p_\varepsilon+q_y)\big)
- f\big(\bar x,\bar t, u(\bar x,\bar t), s_y(p_\varepsilon+q_y)\big),\\
&&  \mathcal{T}_3=
f\big(\bar x,\bar t, u(\bar x,\bar t), s_y(p_\varepsilon+q_y)\big)
-\mu\,  f\big(\bar x,\bar t, u(\bar x,\bar t), \frac{1}{\mu}s_x(p_\varepsilon+q_x)\big),
\end{eqnarray*} 
and we estimate the three terms separately.

From {\bf (B)}(ii), we have
\begin{eqnarray*}
\mathcal{T}_1 \leq m_{R_{\mu,\theta}}\left(
(1+|v(\bar y,\bar t)|+|s_y(p_\varepsilon+q_y)||\bar x-\bar y|
\right),
\end{eqnarray*} 
where $R_{\mu,\theta}$ is given by~\eqref{xybornes}. Using  {\bf (B)}(iv)
and the fact that $v\in \mathcal{C}_{p},$ we get
\begin{eqnarray}\label{ineg-util}
|v(\bar y,\bar t)|, |s_yq_y|=O(R_{\mu,\theta}) 
\quad {\rm and}\quad
|s_y p_\varepsilon||\bar x-\bar y|\leq m(\eps),
\end{eqnarray} 
and therefore
\begin{eqnarray}\label{ineg-654}
\mathcal{T}_1\leq  m(\eps).
\end{eqnarray} 

To deal with $\mathcal{T}_2,$ we first note that
\begin{eqnarray*}
 \mu u(\bar x,\bar t)- v(\bar y,\bar t)
&\geq&
\Phi(\bar x,\bar y, \bar t) \\
&\geq&
 \Phi(\hat x,\hat x, \hat t)\\
&\geq&
\mu u(\hat x,\hat t)-v(\hat y,\hat t)-e^{L\hat t}\theta (1-\mu)(1+2|\hat x|^2)^{p/2}
-\rho \hat t. 
\end{eqnarray*} 
Since $u(\hat x,\hat t)>v(\hat y,\hat t)$ by~\eqref{contr879}, if we take
$\mu$ close enough to 1 and $\rho,\theta$ close enough to 0, we obtain that
$$
\mu u(\bar x,\bar t)\geq v(\bar x,\bar t).
$$ 
From  {\bf (B)}(v) (monotonicity of $f$ in $u$), it follows that 
\begin{eqnarray}
  \mathcal{T}_2&\leq&  f\big(\bar x,\bar t, v(\bar y,\bar t), s_y(p_\varepsilon+q_y)\big)-
   f\big(\bar x,\bar t, \mu u(\bar y,\bar t), s_y(p_\varepsilon+q_y)\big)\nonumber\\
  && + f\big(\bar x,\bar t, \mu u(\bar y,\bar t), s_y(p_\varepsilon+q_y)\big)-
   f\big(\bar x,\bar t,  u(\bar y,\bar t), s_y(p_\varepsilon+q_y)\big)\nonumber\\
   &\le & (1-\mu)|u(\bar x,\bar t)|\nonumber\\
&\leq& C_u(1-\mu)(1+|\bar x|^2)^{p/2},
   \label{ineg-655}
 \end{eqnarray} 
 since $u\in \mathcal{C}_{p}.$

To estimate $\mathcal{T}_3,$ we first 
recall the following convex inequality. If $\Psi :\R^N\to \R$ is
convex and $0<\mu<1,$ then, for all $\xi,\zeta\in\R^N,$ we have
\begin{eqnarray}\label{ineg-convexe}
-\mu\, \Psi (\xi)+ \Psi (\zeta)\leq (1-\mu)\Psi \left(\frac{\mu\xi -\zeta}{\mu-1}\right).
\end{eqnarray}
By {\bf (B)}(iii) (convexity of $f$ with respect to the gradient variable), 
for all $z_1,z_2\in \R^N,$ we obtain
\begin{eqnarray*}
f(x,t,u,z_1)-\mu\, f\left(x,t,u,\frac{z_2}{\mu}\right)
\leq (1-\mu)\, f\left(x,t,u,\frac{z_1-z_2}{1-\mu}\right).
\end{eqnarray*} 
Therefore
\begin{eqnarray}
 \mathcal{T}_3 
&\leq&
(1-\mu)\,  f\left(\bar x,\bar t, u(\bar x,\bar t), \frac{1}{1-\mu}
(s_y(p_\varepsilon+q_y) - s_x(p_\varepsilon+q_x))  \right)\nonumber \\
&\leq&
C_f(1-\mu)\big(1+ |\bar x|^{p}+ |u(\bar x,\bar t)|+ 
\left|\frac{s_y(p_\varepsilon+q_y) - s_x(p_\varepsilon+q_x)}{1-\mu}\right|^{p^\prime}
\big) \label{cc675}
\end{eqnarray} 
by {\bf (B)}(i).
But
\begin{eqnarray*}
s_y(p_\varepsilon+q_y) - s_x(p_\varepsilon+q_x)= (s_y-s_x)p_\varepsilon
+(s_y-s_x)q_y+ s_x(q_y-q_x).
\end{eqnarray*} 
Hence for some $C>0$ depending only on $p$ (which may change during
the computation), we have 
\begin{eqnarray*}
\left|\frac{s_y(p_\varepsilon+q_y) - s_x(p_\varepsilon+q_x)}{1-\mu}\right|^{p^\prime}
&\leq&
\frac{CC_s^{p'}}{(1-\mu)^{p'}}
\left(
(|\bar x-\bar y| |p_\eps|)^{p'} + (|\bar x-\bar y| |q_y|)^{p'}+ |q_x-q_y|^{p'}\right)
\\
&\leq&  e^{p^\prime L\bar t} m(\eps) 
+ CC_s^{p'} e^{p^\prime L\bar t}\theta^{p'}(1+|\bar x|^2+|\bar y|^2)^{p/2},
\end{eqnarray*} 
by using~\eqref{ineg-util}.
Finally, since $u\in \mathcal{C}_{p},$
we get from~\eqref{cc675}
\begin{eqnarray}\label{ineg-656}
 \mathcal{T}_3 
&\leq&
(1-\mu)C_f(1+C_u+CC_s^{p^\prime} e^{p^\prime L\bar t}\theta^{p^\prime})(1+|\bar x|^2+|\bar y|^2)^{p/2}
+ e^{p^\prime L\bar t} m(\eps).
 \end{eqnarray}

\noindent{\it 8. End of the case $\bar t>0,$ choice of the various parameters.}
By plugging
estimates~\eqref{ineg222}, \eqref{ineg-b}, \eqref{ineg-ell},
\eqref{ineg-654}, \eqref{ineg-655} 
and~\eqref{ineg-656} 
in~\eqref{ineg-651}, we get
\begin{eqnarray} 
Le^{L\bar t}\theta(1\!-\!\mu)(1+|\bar x|^2+|\bar y|^2)^{p /2} +\rho
&\leq&
\left( C_1 e^{L\bar t}\theta +  C_2+ C_3 e^{p'L\bar t}\theta^{p'}
\right)(1\! -\! \mu)(1+ |\bar x|^2 + |\bar y|^2)^{p /2}\nonumber \\
&& + \left( -\nu +\theta e^{L\bar t}(C_4 +m(\eps))\right)
(1\!-\!\mu)|\alpha_\eta|^p\nonumber \\
&& + (1+e^{p'L\bar t})m(\eps)+ m(\varrho/\eps^4)+\eta,\label{est-fin1}
\end{eqnarray} 
where
\begin{eqnarray*} 
C_1=4(p - 1)^2C_\sigma^2+4(p + 1)C_b, && C_2=C_\ell  +  C_u +C_f(1+C_u), \\
C_3 = C_fCC_s^{p'}, && C_4= 4(p-1)C_\sigma^2+ C_b,
\end{eqnarray*} 
are positive constants which depend only on the given datas of the problem.

Now we choose the different parameters in order to have a contradiction in the 
above inequality.  We first assume that the
final time $T$ such that
$$
T=1/L>0
$$
(we will recover the result on any interval $[0,T]$ by a step-by-step argument). 
The main difficulty in the above estimate is to deal with the term in 
$|\alpha_\eta|^p$ since the control $\alpha_\eta$ is unbounded. Taking $ \theta >0$ 
such that 
\begin{eqnarray*} 
\theta e^1(C_4+1)\leq \frac{\nu}{2},
\end{eqnarray*} 
we obtain that the coefficient in front of $|\alpha_\eta|^p$ is
negative (we can assume that $\eps$ is small enough in order to have $m(\eps)\leq 1$). 
Then we fix 
\begin{eqnarray} \label{fixe-L}
L> C_1 +\frac{C_2}{\theta}+C_3 e^{p'-1}\theta^{p'-1}
\quad {\rm and} \quad \eta <\frac{\rho}{2}.
\end{eqnarray} 
Therefore~\eqref{est-fin1} implies
\begin{eqnarray*} 
\frac{\rho}{2}\leq (1+e^{p'})m(\eps)+ m(\varrho/\eps^4).
\end{eqnarray*} 
Sending first $\varrho\to 0,$ we obtain a contradiction for small $\eps.$
In conclusion, up to a suitable  choice of the parameters $\theta, L, \eta,$ 
the claim of the Step 3 is proved if $T\leq 1/L.$\\

\noindent{\it 9. Case when $\bar{t}=0.$}
  We have just proved that
there is a subsequence $\eps_n$ such that $\bar t=0$. 
Therefore for $n$ large enough, for all $(x,t)\in\R^N\times [0,T],$
$T\leq 1/L,$ we have 
\begin{eqnarray*} 
 & & 
 \mu u(x,t)-v(x,t)- \theta(1-\mu)e^{Lt}(1+2|x|^2)^{p/2}-\rho t\\ 
 & \le& 
 \mu u(\bar x, 0)-v(\bar y,0)- \theta(1-\mu)(1+|\bar x|^2+|\bar y|^2)^{p/2}-\frac{|\bar x-\bar y|^2}{\eps_n^2}\\ 
 & \le&
 (1-\mu)(|u(\bar x, 0)|- \theta (1+|\bar x|^2)^{p/2})+ u(\bar x, 0)-v(\bar{y},0) \\
 & \le&
(1-\mu)M_\theta+ u(\bar x, 0)-v(\bar{y},0)
 \end{eqnarray*} 
where $M_\theta$ is given by~\eqref{maxC} since $u\in \mathcal{C}_{p}.$
Since $u-v$ is upper-semicontinuous, from~\eqref{def-compxy}, we get
\begin{eqnarray*}
\mathop{\rm lim\,sup}_{\eps_n\to 0}u(\bar x, 0)-v(\bar{y},0) 
\leq u(x_0, 0)-v({x}_0,0)\leq 0,
 \end{eqnarray*} 
using that $u(x_0,0)\leq \psi(x_0)\leq v(x_0,0).$
It follows 
\begin{eqnarray*} 
 \mu u(x,t)-v(x,t)- \theta(1-\mu)e^{Lt}(1+2|x|^2)^{p/2}-\rho t
\leq (1-\mu)M_\theta.
 \end{eqnarray*} 
Sending $\mu\to 1$ and
$\rho\to 0,$ we get $u\le v$ in $\R^N\times [0,T],$ $T\leq 1/L.$
Noticing that $L$ given by~\eqref{fixe-L} depends only on the given constants
of the problem, we recover the comparison on $[0,T]$ for any
$T>0$ by a classical step-by-step argument. 
It completes the proof
of the theorem.
\hfill $\Box$ \\

We end with the proof of the existence result.\\

\noindent{\bf Proof of Theorem~\ref{thm-existence2}.} The point is to
build a sub- and a supersolution. We treat the case of the subsolution
(the case of the supersolution being simpler). It suffices to prove
that, there exists $\tau >0$ such that,
for all $\eps >0$ there exists $M_\eps >0$ such that
\begin{eqnarray}\label{soussol1} 
u_\eps (x,t)= -e^{\rho t}(M_\eps +\eps (1+|x|^p))
\end{eqnarray} 
is a subsolution of~\eqref{hjb} in $\R^N\times [0,\tau]$
with initial data $\psi.$
Indeed, $u_\eps$ does not belong
to $\mathcal{C}_p$ but $\underline u:= {\rm sup}_{\eps >0} u_\eps\in \mathcal{C}_p$
and $\underline u$ is still a subsolution.

Let $\eps >0.$ Since $\psi, \ell(\cdot, t,\alpha), f(\cdot,t,u,z)\in \mathcal{C}_p,$ there exists
$M_\eps =M_\eps(\psi, \ell, f)$ such that $|\psi|, |\chi|, |\gamma|\leq M_\eps+\eps (1+|x|^p).$
Let $u_\eps$ defined by~\eqref{soussol1} with this choice of $M_\eps.$
Let $\alpha\in A.$ In the following computation, $C>0$ is a constant
which depends only on the given datas of the problem and may change
line to line. 
We have, for all $(x,t)\in\R^N\times [0,T],$
\begin{eqnarray*}
\displaystyle \mathcal{L}(u_\eps)&:=&\frac{\partial u_\eps}{\partial t}- 
\langle {b},  Du_\eps\rangle - {\ell} 
-{\rm Trace}\left[ {\sigma}{\sigma}^T 
 D^2u_\eps \right]+f(x,t,u_\eps,sDu_\eps)\\ 
&\leq&
-\rho |u_\eps|
+ C\eps e^{\rho t}(1+|x|+|\alpha|)|x|^{p-1}-\nu|\alpha|^p
+|\chi|+ C\eps e^{\rho t}(1+|x|^2+|\alpha|^2)|x|^{p-2}\\
&& +|\gamma| + C|u_\eps|+ C\eps^{p'} e^{p'\rho t}|x|^{p'(p-1)}\\
&\leq&
-\rho |u_\eps|+ C|u_\eps|+ C\eps^{p'-1} e^{(p'-1)\rho t}|u_\eps|
-\frac{\nu}{2}|\alpha|^p,
\end{eqnarray*}
since $p'(p-1)=p,$
\begin{eqnarray*}
|\alpha||x|^{p-1}+|\alpha|^2|x|^{p-2}\leq
\frac{\nu}{2}|\alpha|^p+C|x|^p
\quad {\rm and} \quad
 |\chi|+|\gamma|\leq 2(M_\eps+\eps(1+|x|^p))=2|u_\eps|.
\end{eqnarray*}
By choosing $\rho$ large enough such that $\rho= C+Ce^1$ and $\tau>0$
such that $(p'-1)\rho\tau\leq 1,$ we obtain $\mathcal{L}(u_\eps)\leq 0$ 
in $\R^N\times [0,\tau].$ Since $u_\eps(\cdot, 0)\leq \psi$ by the
choice of $M_\eps,$ we obtain that $u_\eps$ is a subsolution, which
ends the proof.
\hfill $\Box$

 \section{Equations with superlinear growth on the datas and the solutions}
\label{sec:super}

In this Section we extend the comparison result of~\cite{dll06} for
equations with $p>1$ growth conditions on the datas and on
the solutions. For simplicity, 
we choose to consider here the following model equation where the
diffusion and the drift do  not depend on the control.  

 \begin{equation}\label{model} 
 \left\{\begin{array}{ll} 
 \displaystyle \frac{\partial u}{\partial t}- {\rm Trace}(\sigma\sigma^T D^2 u) +\langle b,Du\rangle
 +f(x,t,u, sDu)=0
&\mbox{in $\R^N\times[0,T],$}\\[2mm]
 u(x,0)=\psi(x) & \mbox{in $\R^N$.} 
  \end{array}\right. 
 \end{equation} 

 The hypothesis on the data are the following:\\[2mm]
 \noindent \textbf{(C)} (Asssumptions on the diffusion and the drift)
\begin{enumerate} 
 
  \item[(i)]~${b}\in C(\R^N \times [0,T]; \R^N)$ and there exists $C_b>0$
such that, for all $x,y\in\R^N,$ $t\in [0,T],$ 
\begin{eqnarray*} 
|{b}(x,t )-{b}(y,t )|& \leq &{C_b} |x-y|, \\ 
|{b}(x,t)| &\leq & {C_b} (1+|x|)\ ; 
\end{eqnarray*} 
 \item[(ii)] 
${\sigma}\in C(\R^N\times [0,T] ; {\cal{M}}_{N,M} (\R))$ is 
Lipschitz continuous with respect to $x$ (uniformly in $t$), 
namely, there exists $C_\sigma >0$ such that, for all $x,y\in \R^N$ and $t\in [0,T],$ 
$$ 
|{\sigma}(x,t)- {\sigma}(y,t)|\leq C_\sigma |x-y|.
$$ 
\end{enumerate} 
Note that $\sigma$ satisfies, { for \ every} $  x\in \R^N, \ t\in [0,T],$ 
\begin{eqnarray*} 
|{\sigma}(x,t)| \leq  {C_{\sigma}}(1+|x| ). 
\end{eqnarray*}

We are able to consider functions which are in a larger class than in
Section~\ref{sec:unbd}.
We say that a function $u\colon\R^N\times[0,T]\to\R$ is in the class $\tilde{\mathcal{C}}_p$ if 
for some $C>0$ we have 
\begin{eqnarray*}\label{cprime}
|u(x,t)|\le C(1+|x|^p),~~~\mbox{for all $(x,t)\in\R^N\times[0,T].$}
 \end{eqnarray*} 

The main result of this Section if the following

 \begin{Theorem}\label{thm-unicbis} 
Assume that $\sigma$ and $b$ satisfy {\bf (C)}, that $f$
satisfies {\bf (B)} and that $\psi\in \tilde{\mathcal{C}}_p.$ 
Let $u \in USC(\R^N\times [0,T])$ be a viscosity subsolution of~\eqref{model} 
and $v \in LSC(\R^N\times [0,T])$ be a viscosity supersolution of~\eqref{model}. 
Suppose that $U$ and $V$ are in the class $\tilde{\mathcal{C}}_p$ and  
satisfy $u(x,0)\leq \psi (x)\leq v(x,0).$  
Then $u\leq v$ in $\R^N\times [0,T].$ 
\end{Theorem} 

Before giving the proof of the theorem, we state an existence result
and provide some examples.
As observed in Section~\ref{sec:unbd}, we can prove the existence of
solutions of~\eqref{hjb} 
(at least for small time) as soon as we are able to build sub- and supersolutions
in the class $\tilde{\mathcal{C}}_p.$ 
In Example \ref{LP-deterministe}, we see that solutions may not exist
for all time.

\begin{Theorem} \label{thm-existence}
Assume \textbf{(B)}--\textbf{(C)}. 
If $K, \rho >0$ are large enough, then $\overline u(x,t)=K{\rm e}^{\rho t}(1+|x|^2)^{p/2}$ is a viscosity
supersolution of~\eqref{model} in $\R^N\times [0,T]$  and
there exists $0< \tau \leq T$ such that
$\underline u(x,t)=-K{\rm e}^{\rho t}(1+|x|^2)^{p/2}$ is a viscosity
subsolution of~\eqref{model} in $\R^N\times [0,\tau].$
In consequence, for all $\psi\in \tilde{\mathcal{C}}_p,$
there exists a unique continuous viscosity solution
of~\eqref{model} in $\R^N\times [0,\tau ]$  in   $\tilde{\mathcal{C}}_p.$
\end{Theorem}

The proof is very close to the one of~\cite[Lemma 2.1]{dll06},
thus  we omit it.
Let us give some examples of Equations for which
Theorem~\ref{thm-unicbis} applies
and some examples.

\begin{Example} \label{exple-simple-p}
The typical (simple) case we have in mind is
\begin{equation}\label{model-p} 
 u_t-\Delta u+|Du|^{p'}= -f(x,t)\quad \mbox{in $\R^N\times[0,T]$},
 \end{equation} 
where $f$ satisfies \textbf{(B)}.
Note that~\eqref{model-p} can be written
$$
u_t-\Delta u+ p\, \mathop{\rm sup}_{\alpha\in\R^N} \{ \langle \alpha , Du \rangle
-\frac{|\alpha|^{p'}}{p'}\} +f(x,t)=0
$$
and therefore is on the form \eqref{hjb}.
The stationary version of this equation was studied in
Alvarez~\cite{alvarez96} under more restrictive assumptions
on the datas and the growth of the solution. More precisely,
he assumed conditions like \eqref{classeC} and \eqref{assumpt-sup}. 
\end{Example}

\begin{Example} \label{expl-bdse}
Equation \eqref{hjb} typically appears in the study of BSDEs where  
$s(x,t)=\sigma(x,t).$ In~\cite{bh08}, Briand and Hu proved that $u$ given 
by~\eqref{nlfk} is a viscosity solution of  \eqref{hjb} for $1\leq
p<2.$ 
Theorem~\ref{thm-unicbis} proves this solution is unique. We are able to deal with any $p>1$ but
we had to impose the regularity condition {\bf (B)}(ii)
on $x$ for $f,$ which is not needed for the BDSEs.
\end{Example}

\begin{Example} \rm 
As far as the coefficient $f$ is concerned, a typical case we have in
mind is
$$
f(x,t,u,z)=g(x,t,u)+|z|^{p'},
$$
with continuous $g$
satisfying {\bf (B)}(i),(ii) and (v).
It leads to nonlinearities
like ``$g(x,t,u)+|s(x,t)Du|^{p'}$'' in the equation.
Note that the power-$p'$ term depends on $x$ via $s(x,t).$
This dependence brings an additional difficulty,
see Lemma~\ref{linearisation}. 
 \end{Example}

\begin{Example} [Deterministic Control Problem]\rm \label{LP-deterministe}
 Consider the control problem (in dimension 1 for sake of simplicity)
\begin{eqnarray*}
\left\{
\begin{array}{cc}
dX_s = {\alpha}_s\, ds, & s\in [t,T], \ 0\leq t\leq T, \\
X_t = x\in \R, &
\end{array}
\right.
\end{eqnarray*}
where the control ${\alpha} \in {\cal{A}}_t := L^p ([t,T] ; \R)$
and the value function is given by
\begin{eqnarray*}
V(x,t)= \mathop{\rm inf}_{{\alpha} \in {\cal{A}}_t}
\{  \int_t^T (\frac{|{\alpha}_s|^p}{p} +\rho |X_s|^p)\, ds - |X_T|^p\} \ \
\ {\rm for \ some} \ \rho >0.
\end{eqnarray*}
 The Hamilton-Jacobi equation formally associated to this problem is 
\begin{eqnarray*} \label{hjb1-lq}
\left\{
\begin{array}{ll}
-  w_t +\frac{1}{p^\prime} | w_x|^{p^\prime}= \rho |x|^{p} \ \
& {\rm in} \ \R\times (0,T), \\
 w (x,T)= -|x|^{p} & {\rm in} \ \R.
\end{array}
\right.
\end{eqnarray*}
 Looking for a solution
$ w$ under the form $ w (x,t)= \varphi (t)|x|^p,$ we obtain that
$\varphi$ is a solution of the differential equation
\begin{eqnarray*}\label{edo-control}
-\varphi ' + \frac{|\varphi|^{p^\prime}}{p^\prime}=\rho \ \  {\rm in} \ (0,T), \ \ \ \ \
\varphi (T)= -1.
\end{eqnarray*}
We get
$$\int_{-1}^{\varphi(t)}\frac{p^\prime}{|y|^{p^\prime}-\rho p^\prime}dy=t-T.$$
One can check that if   $0<\rho p^\prime<1$ 
and $T>\int_{-\infty}^{-1}\frac{p^\prime}{|y|^{p^\prime}-\rho
  p^\prime}dy$,  
then there is $\tau\in (0,T)$ such that the solution
blows up at $t=\tau.$ 
 \end{Example}

Let us turn to the proof of the comparison theorem.\\

\noindent {\bf Proof of Theorem~\ref{thm-unicbis}}. 
To avoid a lot of technicality, we start the proof 
with several lemmas collecting the main 
intermediate results. The proofs of the lemmas are
postponed at the end of the section and can be skipped
at first reading.

\begin{Lemma}\label{chgt-fct} (Change of functions)\\
Let $\tilde{u} ={\rm e}^{-Lt}u+h(x)$ where $h(x)
=\overline{C}(1+|x|^p)$ for some constants
$\overline{C},L >0.$
Then $\tilde{u}$ is a viscosity solution of
\begin{eqnarray}\label{eq-tilde}
\left\{
\begin{array}{ll}
\tilde{u}_t
- {\rm Trace}(\sigma(x,t)\sigma(x,t)^T D^2\tilde{u})+\langle b(x,t),D\tilde{u}\rangle &\\
\hspace*{2cm} +\tilde{f}(x,t,\tilde{u}-h, 
s(x,t)(D\tilde{u}-D{h}))=0
& {in} \ \R^N\times (0,T],\\[2mm]
\tilde{u}(x,0)= \psi(x)+h(x) & {for \ all } \ x\in\R^N,
\end{array}
\right.
\end{eqnarray}
with, for all $(x,t,v,z)\in \R^N\times [0,T]\times \R^N,$
\begin{eqnarray}\label{form-tilde-h}
\label{form-tilde-f}
\tilde{f}(x,t,v,z)=  Lv+\tilde{g}(x,t)
+{\rm e}^{-Lt} f\left(x,t,{\rm e}^{Lt} v, {\rm e}^{Lt}z\right),
\end{eqnarray}
where
\begin{eqnarray}\label{form-tilde-g}
\tilde{g}(x,t)= 
{\rm Trace}(\sigma(x,t)\sigma(x,t)^TD^2h(x))
- \langle b(x,t),Dh(x)\rangle.
\end{eqnarray}
Moreover,
\begin{eqnarray} \label{hyp-tilde-f} 
\tilde{f}(x,t,v,z)
-\tilde{f}(x,t,v',z)\leq (\hat{C}-L)(v'-v)
\quad {if} \ v\leq v'. 
\end{eqnarray} 
\end{Lemma}

In the sequel, since $u,v,\psi\in \tilde{\cal{C}}_p,$
we can choose $\overline{C}>0$ such that
\begin{eqnarray}\label{fixCbarre} 
|u|, |v|, |\psi|\leq \frac{\overline{C}}{2}(1+|x|^p).
\end{eqnarray}
In this case, note that
\begin{eqnarray*} 
\psi(x)+h(x)=\psi(x)+\overline{C}(1+|x|^p)\geq 0
\end{eqnarray*} 
and the initial data is nonnegative in~\eqref{eq-tilde}.

Moreover, we take 
\begin{eqnarray}\label{fixL}
L>\hat{C} \quad {\rm and} \quad
L> 4p(p-1)NC_\sigma^2 +4pC_b+10\hat{C}
\end{eqnarray}
(the constants $C_\sigma$, $C_b$ and $\hat C$ appear in {\bf (B)}).
The first condition ensures that
the right-hand side of~\eqref{hyp-tilde-f} is nonpositive (i.e.
$v\mapsto\tilde{f}(x,t,v,z)$ is nondecreasing). The second
condition appears naturally in the proof of the following lemma.

\begin{Lemma}\label{linearisation} (A kind of linearization
  procedure)\\
Let $\overline{C}, L>0$ be  such that~\eqref{fixCbarre} and~\eqref{fixL}
hold.
Let $0<\mu<1$ and set $\tilde{w}=\mu \tilde{u}-\tilde{v}.$ 
Then $\tilde{w}$ is a viscosity subsolution
of the variational inequality
\begin{eqnarray} \label{var-ineq}
\left\{
\begin{array}{cc}
{\rm min} \left\{ w, {\cal{L}}[w]\right\}\leq 0 & {\rm in}
\ \R^n\times (0,T),\\
w(\cdot ,0)\leq 0 & {\rm in}\ \R^n,
\end{array}\right.
\end{eqnarray}
where 
\begin{eqnarray}\label{op-lin}
{\cal{L}}[w]&:=& \frac{\partial w}{\partial t} -{\rm Trace}
[\sigma(x,t)\sigma^T (x,t)D^2 w]- {C}_b(1+|x|)|Dw|
+\frac{L}{4}(1-\mu)h(x,t) \nonumber\\
&& 
 - (1-\mu) {\rm e}^{-Lt} {f}\left({x},{t},
0, {\rm e}^{Lt} s({x},{t})(
\frac{ Dw}{\mu -1}
- D{h}({x}))\right)
\end{eqnarray}
and ${h}$ is defined in Lemma \ref{chgt-fct}.
\end{Lemma}

\begin{Lemma} \label{prob-parab} (An auxiliary parabolic problem)\\
Consider, for any $R>0,$ the parabolic problem
\begin{equation}\label{equa-parab}
\left\{\begin{array}{ll}
 \varphi_t -r^2 \varphi_{rr} - r
\varphi_{r} =0
 & {  \mbox{in}}\ \
[0,+\infty)\times (0,T],\\
\varphi (r,0)=   {\rm max} \{0,r-R\} & \mbox{in} \ [0,+\infty).
\end{array}\right.
\end{equation}
Then~\eqref{equa-parab} has a unique solution
$\varphi_R \in C([0,+\infty)\times [0,T]) \cap
C^\infty ([0,+\infty)\times (0,T])$ such that,
for all $t\in (0,T],$ $\varphi_R (\cdot,t)$ is positive, nondecreasing
and convex in $[0,+\infty).$ Moreover, for every $(r,t)\in
[0,+\infty)\times (0,T],$
\begin{eqnarray} \label{derivee-phi}
\varphi_R (r,t)\geq {\rm max} \{0,r-R\} , \ \ \
0\leq \frac{\partial \varphi_R}{\partial r} (r,t)\leq {\rm e}^T 
\ \ \ { and} \ \ \
\varphi_R (r,t) \mathop{\longrightarrow}_{R\to +\infty} 0.
\end{eqnarray}
\end{Lemma}
For the proof of Lemma \ref{prob-parab} we refer the reader to~\cite{dll06}.

\begin{Lemma} \label{lem-sursol-stricte}
(Construction of a smooth strict supersolution)\\
Let 
$\Phi (x,t)= \varphi_R (h(x), Ct)$ where
$\varphi_R$ is given by Lemma \ref{prob-parab},
$h(x)= \overline{C}(1+|x|^p),$ $\overline{C}$ satisfies~\eqref{fixCbarre} and $C>0.$ 
Then, for $C$ and $L=L(\mu)$ large enough, we have
\begin{eqnarray}\label{strictsup}
{\cal{L}}[\Phi (x,t)]>0 \quad { for \ all} \ (x,t)\in \R^N\times (0,1/L],
\end{eqnarray}
where ${\cal{L}}$ is defined by~\eqref{op-lin}. 
\end{Lemma}

Now, we continue the {\bf proof of Theorem \ref{thm-unicbis}}.
 Consider
\begin{eqnarray}\label{max-sous-sur}
\mathop{\rm max}_{\R^N\times [0,1/L]} \{ \tilde{w}-\Phi \},
\end{eqnarray}
where $\tilde{w}$ is given by Lemma~\ref{linearisation} and
$\Phi$  is the function  built in Lemma~\ref{lem-sursol-stricte}.
From~\eqref{derivee-phi}, for $|x|$ large enough, we have
$\Phi (x,t) \geq \overline{C}(1+|x|^p)-R.$
Since $\tilde{w}\leq (\mu +1)\overline{C}(1+|x|^p)/2,$
it follows that the maximum~\eqref{max-sous-sur}
is achieved at a point $(\bar{x},\bar{t})\in \R^N\times [0,1/L].$
We can assume that $\tilde{w}(\bar{x},\bar{t})>0$ otherwise, arguing
as in \eqref{abcd12}-\eqref{abcd34}, we prove $\tilde{w}\leq 0$
in $\R^N\times [0,1/L]$ and the conclusion follows.

We claim that  $\bar{t}=0$. Indeed suppose by contradiction that
$\bar{t}>0$. 
Then  since $\tilde{w}$
is a viscosity subsolution of~\eqref{var-ineq} with
$\tilde{w}(\bar{x},\bar{t})>0,$ 
by taking $\Phi$ as a
test-function, we would have
${\cal{L}}[\Phi ](\bar{x},\bar{t}) \leq 0$ which contradicts
the fact that $\Phi$ satisfies (\ref{strictsup})\,.
Thus, for all $(x,t)\in \R^N\times [0,1/L],$
\begin{eqnarray}\label{abcd12}
\tilde{w}(x,t)-\Phi (x,t)\leq \tilde{w}(\bar{x},0)- \Phi (\bar{x},0)
\leq 0,
\end{eqnarray}
where the last inequality follows from~\eqref{var-ineq} and the fact that $\Phi \geq 0.$
Therefore, for every $(x,t)\in \R^N\times [0,1/L],$ we have
\begin{eqnarray}\label{abcd34}
(\mu \tilde{u}-\tilde{v})(x,t) = \tilde{w}(x,t)\leq \Phi(x,t) =\varphi_R(h(x),Ct).
\end{eqnarray}
Letting $R$ to $+\infty,$ we get by~\eqref{derivee-phi},
$\mu \tilde{u}-\tilde{v} \leq 0$ in $\R^N\times [0,1/L].$

We can repeat the above arguments on $\R^N\times [1/L,2/L]$ with
the same constants. By a step-by-step argument, we then prove that
$\mu \tilde{u}-\tilde{v}\leq 0$ in $\R^N\times [0,T].$
Letting $\mu$ go to 1, we obtain $u\leq v$ as well which completes
the proof of the theorem.
\hfill$\Box$\\

We turn to the proof of the Lemmas~\ref{chgt-fct},~\ref{linearisation}
 and~\ref{lem-sursol-stricte}. \\


\noindent {\bf Proof of Lemma~\ref{chgt-fct}}. Since
\begin{eqnarray*}
u= {\rm e}^{Lt}(\tilde{u}-h), 
&\quad&
u_t= {\rm e}^{Lt}(\tilde{u}_t+L(\tilde{u}-h)),\\
Du= {\rm e}^{Lt}(D\tilde{u}-Dh),
&\quad&
D^2u= {\rm e}^{Lt}(D^2\tilde{u}-D^2h),
\end{eqnarray*}
we obtain easily that $\tilde{u}$ is a viscosity solution 
of~\eqref{eq-tilde} with $\tilde{f}$ and $\tilde{g}$ given 
by~\eqref{form-tilde-f} and~\eqref{form-tilde-g}. It remains to 
check~\eqref{hyp-tilde-f}). 
Take $v,v'\in\R$ such that $v\leq v'.$ From 
{\bf (B)}(v), we obtain
\begin{eqnarray*}
&&\tilde{f}(x,t,v,z)
-\tilde{f}(x,t,v',z)\\
&\leq&
L(v-v')+{\rm e}^{-Lt}\left(f(x,t,{\rm e}^{Lt}v,{\rm e}^{Lt}z)
- f(x,t,{\rm e}^{Lt}v',{\rm e}^{Lt}z)\right)\\
&\leq&
-L(v'-v)+{\rm e}^{-Lt}\hat{C}
|{\rm e}^{Lt}(v-v')|
\leq (\hat{C}-L)(v'-v). 
\end{eqnarray*}
It ends the proof of the lemma.
\hfill$\Box$ \\

\noindent{\bf Proof of Lemma~\ref{linearisation}.} 
For $0< \mu <1,$ let $\tilde{u}^{\mu}=\mu\tilde u$ and
$\tilde{w}=\tilde{u}^{\mu}-\tilde{v}.$
We divide the proof in different steps. 

\noindent{\it Step 1. A new equation for $\tilde{u}^{\mu}.$}
It is not difficult to see that,
if $\tilde{u}$ is a subsolution of~\eqref{eq-tilde}, then $\tilde{u}^{\mu}$ is a
subsolution of 
\begin{eqnarray}\label{eq-u-mu}
\left\{
\begin{array}{ll}
\tilde{u}^{\mu}_t
- {\rm Trace}(\sigma(x,t)\sigma(x,t)^T D^2\tilde{u}^{\mu})+\langle b(x,t),D\tilde{u}^{\mu}\rangle &\\
\hspace*{1cm} +\mu\,\tilde{f}(x,t,\displaystyle\frac{\tilde{u}^\mu}{\mu}-{h}, 
s(x,t)(\displaystyle\frac{D\tilde{u}^\mu}{\mu}-D{h}))=0
& {\rm in} \ \R^N\times (0,T],\\[2mm]
\tilde{u}^{\mu}(x,0)= \mu\psi(x)+\mu h(x) & {\rm for \ all } \ x\in\R^N.
\end{array}
\right.
\end{eqnarray}

\noindent{\it Step 2. The equation for $\tilde{w}.$}
Let $\varphi\in C^2 (\R^N\times [0,T])$ and
suppose that we have
\begin{eqnarray}\label{maxw456}
\mathop{\rm max}_{\R^N\times [0,T]} \tilde{w} -\varphi =  (\tilde{w} -\varphi)(\bar{x},\bar{t}).
\end{eqnarray}
We distinguish 3 cases. 

At first,
if the maximum is achieved for $\bar{t}=0,$ then, writing that
$\tilde{u}^\mu$ is a
subsolution of~\eqref{eq-u-mu} and $\tilde{v}$ a supersolution 
of~\eqref{eq-tilde}
at $t=0$ we obtain $\tilde{u}^\mu(\bar{x},0)\leq
\mu\psi(\bar{x})+\mu h(\bar{x})$
and $\tilde{v}(\bar{x},0)\geq
\psi(\bar{x})+h(\bar{x}).$
It follows that
\begin{eqnarray*}
\tilde{w}(\bar{x},0)\leq
(\mu-1)(\psi(\bar{x})+h(\bar{x}))=
(\mu-1)(\psi(\bar{x})+\overline{C}(1+|\bar{x}|^p))\leq 0
\end{eqnarray*}
by~\eqref{fixCbarre}.
Therefore $\tilde{w}$ satisfies~\eqref{var-ineq} at $(\bar{x},0).$

Secondly, we suppose that $\bar{t}>0$ and $\tilde{w}(\bar{x},\bar{t})\leq 0.$
Again,  $\tilde{w}$ satisfies~\eqref{var-ineq} at $(\bar{x},\bar{t}).$

From now on, we consider the last and most difficult case when  
\begin{eqnarray}\label{case-difficile}
\bar{t}>0 \ \ \ {\rm and} \ \ \
\tilde{w}(\bar{x},\bar{t})> 0.
\end{eqnarray}

\noindent{\it Step 3. Viscosity inequalities for $\tilde{u}^\mu$ and $\tilde{v}.$} 
This step is classical in viscosity theory.
We can assume that the maximum in~\eqref{maxw456} at $(\bar{x},\bar{t})$ is strict
in the some ball $\overline{B}(\bar{x},r)\times [\bar{t}-r,\bar{t}+r]$
(see~\cite{barles94} or~\cite{bcd97}). Let
$$
\Theta (x,y,t) = \varphi (x,t)
+\frac{|x-y|^2}{\varepsilon^2}
$$
and consider
$$
M_\varepsilon := \mathop{\rm max}_{x,y\in \overline{B}(\bar{x},r), \,
t\in [\bar{t}-r,\bar{t}+r]}
\{ \tilde{u}^\mu(x,t)- \tilde{v} (y,t)-\Theta (x,y,t) \}.
$$
This maximum is achieved at a point $(x_\varepsilon,
y_\varepsilon, t_\varepsilon)$ and, since the maximum is strict,
we know that
\begin{eqnarray} \label{est-classique3}
x_\varepsilon , y_\varepsilon \to \bar{x}, \quad
\frac{|x_\varepsilon -y_\varepsilon|^2}{\varepsilon^2}\to 0, 
\end{eqnarray}
and
$$
M_\varepsilon = \tilde{u}^\mu(x_\varepsilon,t_\varepsilon) - \tilde{v}
(y_\varepsilon,t_\varepsilon) -\Theta
(x_\varepsilon,y_\varepsilon,t_\varepsilon) \
\mathop{\longrightarrow}  (\tilde{w}
-\varphi) (\bar{x},\bar{t}) \quad
{\rm as} \ \varepsilon\to 0.
$$
It means that, at the limit $\varepsilon\to 0,$ we obtain some
information on $\tilde{w} -\varphi$ at $(\bar{x},\bar{t})$ which will
provide the new equation for $\tilde{w}.$ 
From~\eqref{case-difficile}, for $\eps$ small enough, we have
\begin{eqnarray} \label{strict-pos7}
\tilde{u}^\mu(x_\varepsilon,t_\varepsilon) - \tilde{v}
(y_\varepsilon,t_\varepsilon)>0.
\end{eqnarray}
We can take
$\Theta$ as a test-function to use the fact that $\tilde{u}^\mu$ is a
subsolution of~\eqref{eq-u-mu} and $\tilde{v}$ a supersolution of~\eqref{eq-tilde}.
Indeed $(x,t)\in \overline{B}(\bar{x},r)\times [\bar{t}-r,\bar{t}+r]
\mapsto \tilde{u}^\mu(x,t)- \tilde{v} (y_\varepsilon,t)-\Theta (x,y_\varepsilon,t)$
achieves its maximum at $(x_\varepsilon,t_\varepsilon)$
and $(y,t)\in \overline{B}(\bar{x},r)\times  [\bar{t}-r,\bar{t}+r]
\mapsto -\tilde{u}^\mu (x_\varepsilon,t)+ \tilde{v} (y,t)
+\Theta (x_\varepsilon,y,t)$
achieves its minimum at $(y_\varepsilon,t_\varepsilon).$
Thus, by  Theorem 8.3 in the User's guide~\cite{cil92}, for
every $\rho >0,$ there exist
$a_1, a_2 \in\R$ and $X,Y\in{\mathcal{S}}_{N}$ such that
\begin{eqnarray*}
 \left(a_1,D_x\Theta(x_\varepsilon,y_\varepsilon,t_\varepsilon),
X\right)
  \in\bar{\mathcal{P}}^{2,+}(\tilde{u}^\mu)(x_\varepsilon,
  t_\varepsilon), 
\quad
\left(a_2,-D_y \Theta(x_\varepsilon,y_\varepsilon,t_\varepsilon),  
Y\right)
  \in\bar{\mathcal{P}}^{2,-}( \tilde{v})(y_\varepsilon,t_\varepsilon),
\end{eqnarray*}
$a_1-a_2=\Theta_t (x_\varepsilon,y_\varepsilon,t_\varepsilon)= \varphi_t (x_\varepsilon,t_\varepsilon)$ and
\begin{equation}\label{Ineq-Matrice}
-(\frac{1}{\rho}+|{{M}}|)I \leq
\left(
\begin{array}{cc}
X & 0\\
0 & -Y\\
\end{array}
\right)\leq {{M}}+ \rho {{M}}^2
\ \ \ {\rm where} \
{{M}}=D^2\Theta (x_\varepsilon,y_\varepsilon,t_\varepsilon ).
\end{equation}
Setting $\displaystyle{p_\varepsilon =
2\frac{x_\varepsilon-y_\varepsilon}{\varepsilon^2}},$  we have
\begin{eqnarray*}
D_x\Theta(x_\varepsilon,y_\varepsilon,t_\varepsilon)
=  p_\varepsilon + D \varphi (x_\varepsilon,t_\varepsilon)
\ \ \ {\rm and} \ \ \
 D_y\Theta(x_\varepsilon,y_\varepsilon,t_\varepsilon)
= -p_\varepsilon,
\end{eqnarray*}
and
\begin{eqnarray*}
{{M}}= \left(
\begin{array}{cc}
D^2 \varphi (x_\varepsilon,t_\varepsilon)+2I/\varepsilon^2
& -2I/\varepsilon^2\\
-2I/\varepsilon^2 & 2I/\varepsilon^2\\
\end{array}
\right).
\end{eqnarray*}
Thus, from~\eqref{Ineq-Matrice}, it follows
\begin{eqnarray} \label{ineq-fq}
\langle X p,p\rangle - \langle Y q,q\rangle \leq \langle D^2
\varphi (x_\varepsilon,t_\varepsilon) p,p \rangle
+\frac{2}{\varepsilon^2} |p-q|^2 +
m\left(\frac{\rho}{\varepsilon^4}\right),
\end{eqnarray}
where $m$ is a modulus of continuity which is independent of
$\rho$ and $\varepsilon.$ In the sequel, $m$ will always
denote a generic modulus of continuity independent of
$\rho$ and $\varepsilon.$

Writing the subsolution viscosity inequality for $\tilde{u}^\mu$
and the supersolution inequality for $\tilde{v}$ by
means of the semi-jets and subtracting the inequalities,
we obtain
\begin{eqnarray} \label{eqW}
&& \varphi_t (x_\varepsilon,t_\varepsilon) \nonumber\\
&&- {\rm Trace}\left[ {\sigma} (x_\varepsilon,t_\varepsilon){\sigma}^T(x_\varepsilon,t_\varepsilon)
X \right]+{\rm Trace}\left[ {\sigma} (y_\varepsilon,t_\varepsilon){\sigma}^T(y_\varepsilon,t_\varepsilon)
Y \right]\nonumber\\
&&-
\langle {b}(x_\varepsilon,t_\varepsilon),p_\eps +D\varphi(x_\varepsilon,t_\varepsilon)\rangle+ \langle {b}(y_\varepsilon,t_\varepsilon),p_\eps \rangle \nonumber\\
&&+\mu \tilde{f}\left(x_\eps,t_\eps, 
\frac{\tilde{u}^\mu(x_\eps,t_\eps)}{\mu}-{h}(x_\eps),
s(x_\varepsilon,t_\varepsilon)(\frac{p_\eps+D\varphi(x_\varepsilon,t_\varepsilon)}{\mu}-
D{h}(x_\varepsilon))\right) \nonumber\\
&&- \tilde{f}\left(y_\eps,t_\eps, 
\tilde{v}(y_\eps,t_\eps)-{h}(y_\eps),
s(y_\varepsilon,t_\varepsilon)(p_\eps-
D{h}(y_\varepsilon))\right) \nonumber\\
&\leq& 0
\end{eqnarray}

Now, we derive some estimates for the various terms appearing 
in~\eqref{eqW} in order to be able ton send $\eps\to 0.$ 
The estimates for the $\sigma$ and $b$ terms are
classical wheras those for the $f$ terms are more involved.

For the sake of simplicity, for any function 
$g:\R^N\times [0,T]\to \R,$ we set
$$
g (x_\varepsilon,t_\varepsilon)= {g}_x
\ \ \ {\rm and} \ \ \
{g} (y_\varepsilon,t_\varepsilon)= {g}_y.
$$

\noindent{\it Step 4.  Estimate of $\sigma$-terms.}
Let us denote by $(e_i)_{1\leq i\leq N}$ the canonical basis of $\R^N.$ 
By using~\eqref{ineq-fq}, we obtain
\begin{eqnarray} \label{est-sig}
{\rm Trace}\left[ {c}_x {\sigma}_x^T X -{\sigma}_y {\sigma}_y^T Y \right]
& = & \sum_{i=1}^N
\langle X{\sigma}_x e_i , {\sigma}_x e_i\rangle
-  \langle Y{\sigma}_y e_i , {\sigma}_y e_i\rangle \nonumber\\
& \leq &  
{\rm Trace}\left[ {\sigma}_x {\sigma}_x^T D^2 \varphi (x_\varepsilon,t_\varepsilon )
\right]
+\frac{2}{\varepsilon^2} |{\sigma}_x -{\sigma}_y|^2
+m\left(\frac{\rho}{\varepsilon^4}\right)
\nonumber\\
&\leq &
{\rm Trace}\left[ {\sigma}_x {\sigma}_x^T D^2 \varphi (x_\varepsilon,t_\varepsilon )
\right] +
2C_{{\sigma} ,r}^2 \frac{|x_\varepsilon -y_\varepsilon|^2}{\varepsilon^2}
+ m\left(\frac{\rho}{\varepsilon^4}\right) \nonumber\\
&\leq &
{\rm Trace}\left[ {\sigma} {\sigma}^T(\bar{x},\bar{t}) D^2 \varphi (\bar{x},\bar{t})
\right] +
m(\varepsilon) + m\left(\frac{\rho}{\varepsilon^4}\right),
\end{eqnarray}
where $C_{{\sigma}, r}$ is a Lipschitz constant for ${\sigma}$ in
$\bar{B}(x,r)$ and we used that $\sigma$ is continuous, $\varphi$
is $C^2$ and~\eqref{est-classique3}.

\noindent{\it Step 5.  Estimate of $b$-terms.}
From \textbf{(C)}, if $C_{{b},r}$ is the Lipschitz constant
of ${b}$ in $\overline{B}(\bar{x},r)\times[\bar{t}-r,\bar{t}+r],$
then we have
\begin{eqnarray*}
\langle {b}(x_\varepsilon,t_\varepsilon)
-{b}(y_\varepsilon,t_\varepsilon), p_\varepsilon
\rangle \leq C_{{b},r} |x_\varepsilon -y_\varepsilon| |p_\varepsilon|
\leq 2C_{{b},r} \frac{|x_\varepsilon -y_\varepsilon|^2}{\varepsilon^2}
= m(\varepsilon)
\end{eqnarray*}
and
\begin{eqnarray*}
\langle {b}(x_\varepsilon,t_\varepsilon ),
D\varphi (x_\varepsilon,t_\varepsilon) \rangle
\leq {C}_b(1+|x_\varepsilon|)
|D\varphi (x_\varepsilon,t_\varepsilon)|.
\end{eqnarray*}
It follows
\begin{eqnarray} \label{estb2}
\langle {b}(x_\varepsilon,t_\varepsilon),p_\eps
+D\varphi(x_\varepsilon,t_\varepsilon)\rangle
- \langle {b}(y_\varepsilon,t_\varepsilon),p_\eps \rangle
\leq 
{C}_b(1+|\bar{x}|)|D\varphi (\bar{x},\bar{t})|+
m(\varepsilon)
\end{eqnarray}

\noindent{\it Step 6.  Estimate of $\tilde{f}$-terms.}
We write
\begin{eqnarray*}
&& -\mu \tilde{f}\left(x_\eps,t_\eps, 
\frac{\tilde{u}^\mu(x_\eps,t_\eps)}{\mu}-{h}_x,
s_x(\frac{p_\eps+D\varphi(x_\varepsilon,t_\varepsilon)}{\mu}-
D{h}_x)\right) \\
&& + \tilde{f}\left(y_\eps,t_\eps, 
\tilde{v}(y_\eps,t_\eps)-{h}_y,
 s_y(p_\eps- D{h}_y)\right) \\
&=&  {\cal{T}}_1+{\cal{T}}_2+{\cal{T}}_3
\end{eqnarray*}
where
\begin{eqnarray*}
{\cal{T}}_1 &=&
-\mu \tilde{f}\left(x_\eps,t_\eps, 
\frac{\tilde{u}^\mu(x_\eps,t_\eps)}{\mu}-{h}_x,
s_x(\frac{p_\eps+D\varphi(x_\varepsilon,t_\varepsilon)}{\mu}-
D{h}_x)\right)\\
&&
+ \mu \tilde{f}\left(x_\eps,t_\eps, 
\tilde{v}(y_\eps,t_\eps)-{h}_y,
s_x(\frac{p_\eps+D\varphi(x_\varepsilon,t_\varepsilon)}{\mu}-
D{h}_x)\right),
\\
{\cal{T}}_2 &=&
- \mu \tilde{f}\left(x_\eps,t_\eps, 
\tilde{v}(y_\eps,t_\eps)-{h}_y,
s_x(\frac{p_\eps+D\varphi(x_\varepsilon,t_\varepsilon)}{\mu}-
D{h}_x)\right)
\\
&&
+ \mu \tilde{f}\left(y_\eps,t_\eps, 
\tilde{v}(y_\eps,t_\eps)-{h}_y,
s_x(\frac{p_\eps+D\varphi(x_\varepsilon,t_\varepsilon)}{\mu}-
D{h}_x)\right),
\\
{\cal{T}}_3 &=&
- \mu \tilde{f}\left(y_\eps,t_\eps, 
\tilde{v}(y_\eps,t_\eps)-{h}_y,
s_x(\frac{p_\eps+D\varphi(x_\varepsilon,t_\varepsilon)}{\mu}-
D{h}_x)\right)
\\
&&
+\tilde{f}\left(y_\eps,t_\eps, 
\tilde{v}(y_\eps,t_\eps)-{h}_y,
 s_y(p_\eps-D{h}_y)\right).
\end{eqnarray*}

We estimate ${\cal{T}}_1$.
From~\eqref{strict-pos7}, we have
$$
\tilde{u}(x_\eps,t_\eps)=\frac{\tilde{u}^\mu(x_\eps,t_\eps)}{\mu}
> \tilde{v}(y_\eps,t_\eps) +(1-\mu)\tilde{u}(x_\eps,t_\eps).
$$ 
Using~\eqref{hyp-tilde-f} (the monotonicity  in $u$ of $\tilde{f}$)
and then  {\bf (B)}(v) (Lipschitz continuity in $u$ of $f$), 
we get
\begin{eqnarray*}
&& \tilde{f}\left(x_\eps,t_\eps, 
\frac{\tilde{u}^\mu(x_\eps,t_\eps)}{\mu}-{h}_x,
s_x(\frac{p_\eps+D\varphi(x_\varepsilon,t_\varepsilon)}{\mu}-
D{h}_x)\right)\\
&\geq &
\tilde{f}\left(x_\eps,t_\eps, 
\tilde{v}(y_\eps,t_\eps)+
(1-\mu)\tilde{u}(x_\eps,t_\eps)-{h}_x,
s_x(\frac{p_\eps+D\varphi(x_\varepsilon,t_\varepsilon)}{\mu}-
D{h}_x)\right)\\
&\geq &
\tilde{f}\left(x_\eps,t_\eps, 
\tilde{v}(y_\eps,t_\eps)+{h}_y,
s_x(\frac{p_\eps+D\varphi(x_\varepsilon,t_\varepsilon)}{\mu}-
D{h}_x)\right)\\
&& -\hat{C}(1-\mu)|\tilde{u}(x_\eps,t_\eps)|-\hat{C}|{h}_x-{h}_y|.
\end{eqnarray*}
Since ${h}$ is continuous, we have 
$\hat{C}|{h}_x-{h}_y|=m(\eps).$ By~\eqref{fixCbarre} and
since $x_\eps\to \bar{x},$ we obtain
\begin{eqnarray*}
\hat{C}(1-\mu)|\tilde{u}(x_\eps,t_\eps)|\leq
\hat{C}\overline{C}(1-\mu)(1+|x_\eps|^p)=
\hat{C}(1-\mu)h(\bar{x})+m(\eps).
\end{eqnarray*}
Therefore
\begin{eqnarray}\label{estT1}
{\cal{T}}_1\leq \hat{C}(1-\mu)h(\bar{x})+m(\eps).
\end{eqnarray}

The estimate of ${\cal{T}}_2$ relies on {\bf (B)}(ii). 
Setting $Q_\eps={\rm e}^{Lt_\eps} s_x(\frac{p_\eps+D\varphi(x_\varepsilon,t_\varepsilon)}{\mu}-
D{h}_x)$ and recalling that $r$ is defined at the beginning of Step 3,
we have
\begin{eqnarray}
|{\cal{T}}_2 |
&\leq &
\mu |g(x_\eps,t_\eps)-g(y_\eps,t_\eps)| \nonumber \\
&& + \mu {\rm e}^{-Lt_\eps} 
|f(x_\eps,t_\eps, {\rm e}^{Lt_\eps} (\tilde{v}(y_\eps,t_\eps)-{h}_y),Q_\eps)
- f(y_\eps,t_\eps, {\rm e}^{-Lt_\eps} (\tilde{v}(y_\eps,t_\eps)-{h}_y),Q_\eps)| \nonumber \\
&\leq &
\mu |g(x_\eps,t_\eps)-g(y_\eps,t_\eps)| \nonumber \\
&&+ \mu {\rm e}^{-Lt_\eps}\, m_{2r} \left( (1
+|{\rm e}^{Lt_\eps}(\tilde{v}(y_\eps,t_\eps)-{h}_y)|
+ |Q_\eps| ) 
|x_\eps-y_\eps|\right) \nonumber \\
&\leq & m(\eps),\label{estT2}
\end{eqnarray}
since $g$ is continuous, $|x_\eps-y_\eps|=m(\eps)$ and
$p_\eps|x_\eps-y_\eps|=|x_\eps-y_\eps|^2/\eps^2=m(\eps)$ 
by~\eqref{est-classique3}.

Let us turn to the estimate of ${\cal{T}}_3.$
We have
\begin{eqnarray*}
{\cal{T}}_3 &=&
L (1-\mu) (\tilde{v}(y_\eps,t_\eps)-{h}_y) + (1-\mu)\tilde{g}(y_\eps,t_\eps)\\
&& 
- \mu {\rm e}^{-Lt_\eps}{f}\left(y_\eps,t_\eps, 
{\rm e}^{Lt_\eps}(\tilde{v}(y_\eps,t_\eps)-{h}_y),
{\rm e}^{Lt_\eps}s_x(\frac{p_\eps+D\varphi(x_\varepsilon,t_\varepsilon)}{\mu}-
D{h}_x)\right)
\\
&&
+{\rm e}^{-Lt_\eps} {f}\left(y_\eps,t_\eps, 
{\rm e}^{Lt_\eps}(\tilde{v}(y_\eps,t_\eps)-{h}_y),
 {\rm e}^{Lt_\eps}s_y(p_\eps-D{h}_y)\right).
\end{eqnarray*}
At first, 
from~\eqref{fixCbarre}, we have
 \begin{eqnarray*}
L (1-\mu ) (\tilde{v}(y_\eps,t_\eps)-{h}_y) 
\leq 
-\frac{L(1-\mu )}{2}h(\bar{x})+m(\eps).
\end{eqnarray*}
Using {\bf (C)} (see~\eqref{form-deriv23} for the details),
a straightforward computation gives an estimate for the continuous
function $\tilde{g}$:
 \begin{eqnarray*}
|\tilde{g}(y_\eps,t_\eps)|\leq (p(p-1)NC_\sigma^2 + pC_b)h(\bar{x})+m(\eps).
\end{eqnarray*}
Now we estimate the $f$-terms.
From  {\bf (B)}(iii) (convexity of ${f}$ with respect
to the gradient variable), we can apply~\eqref{ineg-convexe} to
obtain
\begin{eqnarray*}
&& - \mu {\rm e}^{-Lt_\eps}{f}\left(y_\eps,t_\eps, 
{\rm e}^{Lt_\eps}(\tilde{v}(y_\eps,t_\eps)-{h}_y),
{\rm e}^{Lt_\eps}s_x(\frac{p_\eps+D\varphi(x_\varepsilon,t_\varepsilon)}{\mu}-
D{h}_x)\right) \\
&&
+{\rm e}^{-Lt_\eps} {f}\left(y_\eps,t_\eps, 
{\rm e}^{Lt_\eps}(\tilde{v}(y_\eps,t_\eps)-{h}_y),
 {\rm e}^{Lt_\eps}s_y(p_\eps-D{h}_y)\right)\\
&\leq&
(1-\mu) {\rm e}^{-Lt_\eps} {f}\left(y_\eps,t_\eps, 
{\rm e}^{Lt_\eps}(\tilde{v}(y_\eps,t_\eps)-{h}_y),
Q_\eps \right),
\end{eqnarray*}
where
\begin{eqnarray*}
Q_\eps & =& \frac{{\rm e}^{Lt_\eps}}{\mu -1}\left( (s_x-s_y)p_\eps 
+ s_x D\varphi(x_\varepsilon,t_\varepsilon)
+ s_y D{h}_y- \mu s_x D{h}_x  
\right)\\
&=&
{\rm e}^{L\bar{t}} s(\bar{x},\bar{t})\left(
\frac{ D\varphi(\bar{x},\bar{t})}{\mu -1}
-  D{h}(\bar{x})\right)
 +m(\eps).
\end{eqnarray*}
from {\bf (B)}(iv) and~\eqref{est-classique3}.
From {\bf (B)}(v),~\eqref{fixCbarre} and the continuity of $f,$
it follows
\begin{eqnarray*}
&&  {\rm e}^{-Lt_\eps} {f}\left(y_\eps,t_\eps, 
{\rm e}^{Lt_\eps}(\tilde{v}(y_\eps,t_\eps)-{h}_y), Q_\eps\right)\\
&\leq & 
 {\rm e}^{-L\bar{t}} f\left(\bar{x}, \bar{t}, 0, {\rm e}^{L\bar{t}} s(\bar{x},\bar{t})\left(
\frac{ D\varphi(\bar{x},\bar{t})}{\mu -1}
-  D{h}(\bar{x})\right)\right)+ \frac{3\hat{C}}{2}h(\bar{x}) + m(\eps).
\end{eqnarray*}
Finally,  we obtain
\begin{eqnarray} \label{estT3}
 {\cal{T}}_3 
&\leq &  (1-\mu) \left(-\frac{L}{2}
+ p(p-1)NC_\sigma^2 + pC_b+ \frac{3\hat{C}}{2}
\right) h(\bar{x}) \nonumber \\
&& +
 {\rm e}^{-L\bar{t}} f\left(\bar{x}, \bar{t}, 0, {\rm e}^{L\bar{t}} s(\bar{x},\bar{t})\left(
\frac{ D\varphi(\bar{x},\bar{t})}{\mu -1}
-  D{h}(\bar{x})\right)\right)
+m(\eps).
\end{eqnarray}

\noindent{\it Step 7. End of the proof.}
Combining~\eqref{eqW}, \eqref{est-sig}, \eqref{estb2}, \eqref{estT1},
\eqref{estT2} and \eqref{estT3},
setting $L> 4p(p-1)NC_\sigma^2 +4pC_b+10\hat{C}$
and sending $\rho\to 0$ and then $\eps\to 0,$ we get
\begin{eqnarray*}
&& \varphi_t (\bar{x},\bar{t})
-{\rm Trace}\left[ {\sigma} {\sigma}^T(\bar{x},\bar{t}) D^2 \varphi (\bar{x},\bar{t})
\right]
-C_b (1+|\bar{x}|)|D\varphi (\bar{x},\bar{t})|
+\frac{L}{4}(1-\mu)h(\bar{x})
\\
&& - (1-\mu) {\rm e}^{-L\bar{t}}{f}\left(\bar{x},\bar{t},0,
{\rm e}^{L\bar{t}} s(\bar{x},\bar{t})\left(
\frac{ D\varphi(\bar{x},\bar{t})}{\mu -1}
-  D{h}(\bar{x})\right)\right)
\\
&\leq & 0,
\end{eqnarray*}
which is exactly the new equation for $\tilde{w}$ in the 
case~\eqref{case-difficile}. It completes the proof of the lemma.
~~~\hfill$\Box$\\

\noindent{\bf Proof of Lemma \ref{lem-sursol-stricte}.} 
For simplicity, we fix $R$ and
set $\varphi=\varphi_R$ for simplicity. Therefore
$\varphi_r$ denotes the derivative of $\varphi$ wrt
the space variable. We compute
\begin{eqnarray*}
&&\Phi_t = C\varphi_t,
\quad
D\Phi = \varphi_r Dh, \\
&&
D^2\Phi = \varphi_r  D^2h + \varphi_{rr} Dh\otimes Dh,
\end{eqnarray*}
with
\begin{eqnarray*}
h=\overline{C} (1+|x|^p),
\quad
Dh=p\overline{C}|x|^{p-2}x,
\quad
D^2h= p \overline{C}
(|x|^{p-2} Id + (p-2)|x|^{p-4}x\otimes x).
\end{eqnarray*}
For all $(x,t)\in\R^N\times (0,T],$
\begin{eqnarray}
&&\mathcal{L}(\Phi(x,t))\nonumber \\
&=&
C\varphi_t 
- \left({\rm Trace}(\sigma\sigma^TD^2h)+{C}_b(1+|x|)|Dh|\right)\varphi_r
- {\rm Trace}(\sigma\sigma^TDh\otimes Dh)\varphi_{rr}\nonumber \\
&& +(1-\mu )\frac{L}{4}h
 -(1-\mu ) {\rm e}^{-Lt}{f}\left( x,t,0,
{\rm e}^{Lt}s\left(\frac{\varphi_r}{\mu -1}+1\right)Dh\right).
\label{minor-phi}
\end{eqnarray}
Using {\bf (C)}(ii) and the fact that $p'(p-1)=p,$
we have the following estimates:
\begin{eqnarray}\label{form-deriv23}
\begin{array}{l}
 |Dh|\leq p\overline{C}|x|^{p-1}, \quad
|D^2h|\leq p(p-1)\overline{C}|x|^{p-2},\\[1mm]
{C}_b(1+|x|)|Dh|\leq pC_b h,\\[1mm]
 |{\rm Trace}(\sigma\sigma^TD^2h)|
\leq p(p-1)NC_\sigma^2 h,\\
0\leq {\rm Trace}(\sigma\sigma^TDh\otimes Dh)
\leq C_\sigma^2(1+|x|^{2})p^2\overline{C}^2|x|^{2(p-1)}
\leq p^2C_\sigma^2h^2.
\end{array}
\end{eqnarray}
Now, the assumption {\bf (B)}(i) on the growth of $f$ plays a crucial role:
\begin{eqnarray*}
&& {f}\left( x,t,0,
{\rm e}^{Lt}s\left(\frac{\varphi_r}{\mu -1}+1\right)Dh\right)\\
&\leq &
C_f \left(1+|x|^p + \left|{\rm e}^{Lt}s\left(\frac{\varphi_r}{\mu -1}+1\right)Dh \right|^{p'}
\right)\\
&\leq &
C_f\left( \frac{1}{\overline{C}}+
p^{p'}C_s^{p'}\overline{C}^{p'-1}{\rm e}^{Lp't}
\left(\frac{{\rm e}^T}{1-\mu}+1\right)^{p'} \right)h(x)
\end{eqnarray*}
since $\varphi_r\leq {\rm e}^T$ (Lemma \ref{prob-parab}) and  
$|Dh|^{p'}\leq p^{p'}\overline{C}^{p'-1}h(x)$ (because $p'(p-1)=p$).
It follows from~\eqref{minor-phi},
\begin{eqnarray*}
&&\mathcal{L}(\Phi(x,t))\\
&\geq &
C\varphi_t -( p(p-1)NC_\sigma^2+pC_b)\varphi_r - p^2C_\sigma^2h^2\varphi_{rr} \\
&& +(1-\mu)\left( \frac{L}{4}- \frac{C_f {\rm e}^{-Lt}}{\overline{C}}
- p^{p'}C_s^{p'}\overline{C}^{p'-1}{\rm e}^{Lp't}
\left(\frac{{\rm e}^T}{1-\mu}+1\right)^{p'} \right) h(x).
\end{eqnarray*}
We take
\begin{eqnarray*}
C &>& {\rm max}\left\{ p(p-1)NC_\sigma^2+pC_b, p^2C_\sigma^2\right\},\\
L &>& \frac{4C_f}{\overline{C}}+
4 p^{p'}C_s^{p'}\overline{C}^{p'-1}{\rm e}^{p'}
\left(\frac{{\rm e}^T}{1-\mu}+1\right)^{p'}+1
\quad {\rm and \ \eqref{fixL} \ holds},\\
\tau &=& \frac{1}{L}.
\end{eqnarray*}
For this choice of parameters, for all $(x,t)\in \R^N\times (0,\tau],$
  we have
\begin{eqnarray*}
\mathcal{L}(\Phi(x,t))\geq 
C\left( \varphi_t(h,Ct)- h\varphi_r(h,Ct)- h^2 \varphi_{rr}(h,Ct)
\right)+(1-\mu)h >0
\end{eqnarray*}
since $\varphi$ is a solution of~\eqref{equa-parab} and $h>0.$
This proves the lemma.
\hfill$\Box$\\




 \end{document}